\documentclass{naw-latex}
\usepackage[dutch]{babel}

\usepackage{amsmath}
\usepackage{amssymb}
\usepackage{enumerate}
\usepackage{varioref}
\usepackage{charter,eulervm}
\usepackage{boxedminipage}

\newcommand{\es}{\varnothing}

\newtheorem{theorem}{Stelling}
\newtheorem{definition}{Definitie}
\newtheorem{remark}{Opmerking}
\newtheorem{corollary}{Gevolg}
\newtheorem{lemma}{Lemma} 

\newenvironment{bewijs}%
     {\par\ifdim\lastskip<\smallskipamount \removelastskip\smallskip\fi
      \noindent\textsc{Bewijs}.\enspace\ignorespaces}%
     {\par\smallskip}

\bibitem{kn:aardenne}{van~Aardenne-Ehrenfest,~T. and N.~G.~de~Bruijn, 
Circuits and trees in oriented linear graphs, 
{\em Simon Stevin\/} {\bf 28} (1951), pp.~203--217.\\
Dit artikel is ook verschenen in: (I.~Gessel and G.~Rota eds.) 
{\em Classic Papers in Combinatorics\/}, Modern Birkh\"auser Classics, 
2009, pp.~149--163}

\bibitem{kn:amin}{Amin,~K., Factorization of finite abelian groups, 
{\em International Journal of Algebra\/} {\bf 6} (2012), pp.~101--107.} 

\bibitem{kn:andre}{Andr\'e,~D., D\'eveloppements de $\sec x$ et de 
$\tan x$, 
{\em Comtes Rendus des S\'eances de l'Acad\'emie 
des Sciences. Paris.\/} 88 (1879), pp.~965--967.} 

\bibitem{kn:andre2}{Andr\'e,~D., Sur les permutations altern\'ees, 
{\em Journal de Math\'ematiques pures et appliqu\'ees. Paris.\/} 
{\bf 7} (1881), 
pp.~167--184.} 

\bibitem{kn:andre3}{Andr\'e,~D., Probabilit\'e pour qu'une permutation 
donn\'ee de $n$ lettres soit une permutation altern\'ee, 
{\em Comtes Rendus des S\'eances de l'Acad\'emie 
des Sciences. Paris.\/} {97} (1883), pp.~983--984.} 

\bibitem{kn:berge}{Berge,~C., 
Two theorems in graph theory, 
{\em Proceedings of the National Academy of Sciences of the 
United States of America\/} {\bf 43} (1957), pp.~842--844.} 
   
\bibitem{kn:bjorklund}{Bj\"orklund,~A., 
Counting perfect matchings as fast as Ryser, 
{\em Proceedings $23^{\mathrm{rd}}$ Annual ACM-SIAM Symposium on Discrete 
Algorithms\/}, Kyoto, Japan. SIAM (2012), pp.~914--921.} 

\bibitem{kn:brightwell}{Brightwell,~G. and P.~Winkler, 
Counting Eulerian cirsuits is $\#$P-complete, 
{\em Proceedings of the $7^{\mathrm{th}}$ Workshop on Algorithm Engineering 
and Experiments and the Second Workshop on Analytic Algorithmics and 
Combinatorics\/}, Vancouver, BC, Canada. SIAM (2005), pp.~259--262.}  

\bibitem{kn:bruijn24}{de~Bruijn,~N.~G., 
Gemeenschappelijke representatensystemen van twee klassenindeelingen 
van een verzameling, 
{\em Nieuw Archief voor Wiskunde\/} {\bf 22} (1943), pp.~48--52.} 

\bibitem{kn:bruijn1}{de~Bruijn,~N.~G.,  
A combinatorial problem, 
{\em Idagationes Mathematicae\/} {\bf 8} (1946), pp.~461--467.}

\bibitem{kn:bruijn20}{de~Bruijn,~N.~G., 
On bases for the set of integers, 
{\em Publicationes Mathematicae\/} {\bf 1} (1950), 
pp.~232--242.} 

\bibitem{kn:bruijn5}{de~Bruijn,~N.~G., 
On the factorization of finite abelian groups, 
{\em Indagationes Mathematicae\/} {\bf 15} (1953), pp.~258--264.} 

\bibitem{kn:bruijn6}{de~Bruijn,~N.~G., 
On the factorization of cyclic groups, 
{\em Indagationes Mathematicae\/} {\bf 15} (1953), pp.~370--377.} 

\bibitem{kn:bruijn22}{de~Bruijn,~N.~G., 
On number systems, 
{\em Nieuw Archief voor Wiskunde\/} {\bf 4} (1956), pp.~15--17.} 

\bibitem{kn:bruijn9}{de~Bruijn,~N.~G., 
Generalization of P\'olya's fundamental theorem in enumerative 
combinatorial analysis, 
{\em Indagationes Mathematicae\/} {\bf 21} (1959), pp.~59--69.} 

\bibitem{kn:bruijn26}{de~Bruijn,~N.~G., 
Enumerative combinatorial problems concerning structures, 
{\em Nieuw Archief voor Wiskunde\/} {\bf 11} (1963), pp.~142--161.} 

\bibitem{kn:bruijn21}{de~Bruijn,~N.~G., 
Some direct decompositions of the set of integers, 
{\em Mathematics of Computation\/} {\bf 18} (1964), pp.~537--546.} 
 
\bibitem{kn:bruijn8}{de~Bruijn,~N.~G., 
P\'olya's theory of counting. 
Chapter 5 in (E.~Beckenbach ed.) 
{\em Applied Combinatorial mathematics\/}, Wiley, 1964, pp.~144--184.} 

\bibitem{kn:bruijn12}{de~Bruijn,~N.~G., 
Colour patterns that are invariant under a given permutation of the colours, 
{\em Journal of Combinatorial Theory\/} {\bf 2} (1967), pp.~418--421.} 

\bibitem{kn:bruijn25}{de~Bruijn,~N.~G., 
Filling boxes with bricks, 
{\em American Mathematical Monthly\/} {\bf 76} (1969), pp.~37--40.} 

\bibitem{kn:bruijn23}{de~Bruijn,~N.~G., 
Permutations with given ups and downs, 
{\em Nieuw Archief voor Wiskunde\/} {\bf 18} (1970), pp.~61--65.} 

\bibitem{kn:bruijn10}{de~Bruijn,~N.~G., 
P\'olya's Abz\"ahltheorie: Muster f\"ur Graphen und chemische 
Verbindungen. In (K.~Jacobs ed.) 
{\em Selecta Mathematica III\/} Springer-Verlag, Heidelberger 
Taschenb\"ucher {\bf 86}, 1971, pp.~1--26.}

\bibitem{kn:bruijn11}{de~Bruijn,~N.~G., 
A survey of generalizations of P\'olya's 
enumeration theorem, 
{\em Nieuw Archief voor Wiskunde\/} {\bf 19} (1971), pp.~89--112.} 

\bibitem{kn:bruijn14}{de~Bruijn,~N.~G., 
Enumeration of mapping patterns, 
{\em Journal of Combinatorial Theory (A)\/} {\bf 12} (1972), 
pp.~14--20.} 

\bibitem{kn:bruijn27}{de~Bruijn,~N.~G., 
A solitaire game and its relation to a finite field, 
{\em Journal of Recreational Mathematics\/} {\bf 5} (1972), pp.~133--137.}

\bibitem{kn:bruijn2}{de~Bruijn,~N.~G., 
Acknowledgement of priority to C. Flye Sainte-Marie on the counting 
or circular arrangements of $2^n$ zeros and ones that show 
each $n$-letter word eaxactly once. 
Technical Report, TUE, T.H.-Report 75-WSK-06, Technological University 
Eindhoven, 1975.}   

\bibitem{kn:bruijn13}{de~Bruijn,~N.~G., 
A note on the Cauchy-Frobenius lemma, 
{\em Indagationes Mathematicae\/} {\bf 41} (1979), pp.~225--228.} 

\bibitem{kn:bruijn15}{de~Bruijn,~N.~G., 
De stelling van P\'olya, met toepassing op het tellen van bomen en 
boomvormige molekulen. 
In {\em Vertelling over tellen\/}, Vakantiecursus 34/80, 
Mathematisch Centrum, Amsterdam, 1980.} 

\bibitem{kn:bruijn19}{de~Bruijn,~N.~G., 
Counting complete matchings without using Pfaffians, 
{\em Indagationes Mathematicae\/} {\bf 42} (1980), pp.~361--366.} 

\bibitem{kn:bruijn17}{de~Bruijn,~N.~G., 
Algebraic theory of Penrose's non-periodic 
tilings of the plane (I and II), 
{\em Indiagationes Mathematicae\/} {\bf 43} (1981), pp.~38--66.} 
 
\bibitem{kn:bruijn16}{de~Bruijn,~N.~G., 
Mijn liefste boek: P\'olya-Szeg\"o. 
In {\em Uitgelezen Gezelschap\/}, Equiliber nr. 6, 
Bibliotheek TUE, 1994, pp.~4--8.}  

\bibitem{kn:bruijn3}{de~Bruijn,~N.~G. and P.~Erd\"os, 
On a combinatorial problem, 
{\em Indagationes Mathematicae\/} {\bf 10} (1948), pp.~421--423.\\
In de originele titel staat een tikfout; `combinatioral.'} 

\bibitem{kn:bruijn4}{de~Bruijn,~N.~G. and P.~Erd\"os, 
A colour problem for infinite graphs and a problem 
in the theory of relations, 
{\em Indagationes Mathematicae\/} {\bf 13} (1951), pp.~371--373.} 

\bibitem{kn:bruijn18}{de~Bruijn,~N.~G., D.~E.~Knuth and 
S.~O.~Rice, The average height of planted plane trees. 
In (R.~Read ed.) {\em Graph Theory and Computing\/}, Academic Press, 
1972, pp.~15--22.} 

\bibitem{kn:bruijn7}{de~Bruijn,~N.~G. and B.~J.~M.~Morselt, 
A note on plane trees, 
{\em Journal of Combinatorial Theory\/} {\bf 2} (1967), 
pp.~27--34.} 

\bibitem{kn:compeau}{Compeau,~P.~E.~C., P.~A.~Pevzner and G.~Tesler, 
How to apply de Bruijn graphs to genome assembly, 
{\em Nature Biotechnology\/} {\bf 29} (2011), pp.~987--991.} 

\bibitem{kn:dinitz}{Dinitz,~M., 
Full rank tilings of finite abelian groups, 
{\em SIAM Journal on Discrete Mathematics\/} {\bf 20} (2006), pp.~160--170.} 

\bibitem{kn:dowell}{Dowell,~J., 
Periodic basic sequences, 
{\em Nieuw Archief voor Wiskunde\/} {\bf 16} (1968), pp.~112--115.} 

\bibitem{kn:eigen}{Eigen,~S.~J., Y.~Ito and V.~S.~Prasad, 
Universally bad integers and the 2-adics, 
{\em Journal of Number Theory\/} {\bf 107} (2004), pp.~322--334.} 

\bibitem{kn:entringer}{Entringer,~R.~C., 
A combinatorial interpretation of the Euler and Bernoulli 
numbers, 
{\em Nieuw Archief voor Wiskunde\/} {\bf 14} (1966), pp.~241--246.} 

\bibitem{kn:everett}{Everett,~C.~J. and G.~Whaples, 
Representations of sequences of sets, 
{\em American Journal of Mathematics\/} {\bf 71} (1949), pp.~287--293.} 

\bibitem{kn:flye}{Flye Sainte-Marie,~C., 
Solution to question nr. 48, 
{\em L'Interm\'ediare des Math\'ematiciens\/} {\bf 1} (1894), 
pp.~107--110.}

\bibitem{kn:gallai}{Gallai (Gr\"unwald),~T., 
Solution to problem 4065, 
{\em American Mathematical Monthly\/} {\bf 51} (1944), pp.~169--171.}

\bibitem{kn:good}{Good,~I.~J., 
Normal recurring decimals, 
{\em Journal of the London Mathematical Society\/} {\bf 21} (1946), 
pp.~167--169.} 

\bibitem{kn:gottschalk}{Gottschalk,~W.~H., 
Choice functions and Tychonoff's theorem, 
{\em Proceedings of the American Mathematical Society\/} {\bf 2} 
(1951), pp.~172.} 

\bibitem{kn:hajos}{Haj\'os,~G., 
\"Uber einfache und mehrfache Bedeckung des $n$-dimensionalen 
Raumes met einem W\"urfelgitter, 
{\em Mathematische Zeitschrift\/} {\bf 47} (1941), pp.~427--467.} 

\bibitem{kn:hall}{Hall,~M., 
Distinct representatives of subsets, 
{\em Bulletin of the American Mathematical Society\/} {\bf 54} (1948), 
pp.~922--926.} 

\bibitem{kn:halmos}{Halmos,~P.~R. and H.~E.~Vaughan, 
The marriage problem, 
{\em American Journal of Mathematics\/} {\bf 72} (1950), pp.~214--215.} 
 
\bibitem{kn:harary}{Harary,~F., G.~Prins and W.~T.~Tutte, 
The number of plane trees, 
{\em Indagationes Mathematicae\/} {\bf 26} (1964), pp.~319--329.}

\bibitem{kn:kasteleyn2}{Kasteleyn,~P.~W., 
The statistics of dimers on a lattice. I. The number of 
dimer arrangements on a quadratic lattice, 
{\em Physica\/} {\bf 27} (1961), pp.~1209--1225.} 

\bibitem{kn:kasteleyn3}{Kasteleyn,~P.~W., 
Dimer statistics and phase transitions, 
{\em Journal of Mathematical Physics\/} {\bf 4} (1963), pp.~287--293.} 

\bibitem{kn:kasteleyn}{Kasteleyn,~P.~W., 
Graph theory and crystal physics. 
In (F.~Harary ed.) {\em Graph Theory and Theoretical Physics\/}, 
Academic Press, 1967, pp.~43--110.}  

\bibitem{kn:konig}{K\"onig,~D. and S.~Valk\'o, 
\"Uber mehrdeutige Abbildungen von Mengen, 
{\em Mathematische Annalen\/} {\bf 95} (1926), pp.~135--138.} 

\bibitem{kn:konig2}{K\"onig,~D., 
\"Uber Graphen und ihre Anwendung auf Determinantentheorie und 
Mengenlehre, 
{\em Mathematische Annalen\/} {\bf 77} (1916), pp.~453--465.} 

\bibitem{kn:lint}{van~Lint,~J.~H., 
Recente ontwikkelingen in de combinatoriek, 
{\em Nieuw Archief voor Wiskunde\/} {\bf 24} (1976), pp.~215--225.}  

\bibitem{kn:little}{Little, C.~H.~C.,
Kasteleyn's theorem and arbitrary graphs, 
{\em Canadian Journal of Mathematics\/} {\bf 25} (1973), pp.~758--764.} 

\bibitem{kn:little2}{Little,~C.~H.~C., 
Extensions of Kasteleyn's method of enumerating the 1-factors 
of planar graphs. In (D.~Holton ed.) {\em Combinatorial Mathematics II\/}, 
Proceedings of the Second Australian Conference, Melbourne, Australia.
Springer-Verlag, 
Lecture Notes in Mathematics {\bf 403}, (1974), pp.~63--72.}  
 
\bibitem{kn:little3}{Little,~C.~H.~C., 
A characterization of convertible $(0,1)$-matrices, 
{\em Journal of Combinatorial Theory, Series B\/} {\bf 18} (1975), 
pp.~187--208.} 

\bibitem{kn:macmahon}{MacMahon,~P.~A., 
Second memoir on the compositions of numbers, 
{\em Philosophical Transactions of the Royal Society of London. 
Series A, Containing Papers of a Mathematical or Physical Character\/} 
{\bf 207} (1908), pp.~65--134.}  

\bibitem{kn:macmahon2}{MacMahon,~P.~A., 
{\em Combinatory Analysis\/}, Volume 1,  
Cambridge University Press, 1915.}  
 
\bibitem{kn:moser}{Moser,~L., 
An application of generating series, 
{\em Mathematical Magazine\/} {\bf 35} (1962), pp.~37--38.} 

\bibitem{kn:motzkin}{Motzkin,~Th., 
The lines and planes connecting the points of a finite set, 
{\em Transactions of the American Mathematical Society\/} {\bf 70} 
(1951), pp.~451--464.} 

\bibitem{kn:nienhuys}{Nienhuys,~J.~W., 
(L.~Hung and T.~Kloks, eds.) 
{\em De~Bruijn's Combinatorics\/}. Manuscript on 
viXra: 1208.0223, 2012.}  

\bibitem{kn:niven}{Niven,~I., 
A combinatorial problem on finite sequences, 
{\em Nieuw Archief voor Wiskunde\/} {\bf 16} (1968), pp.~116--123.} 

\bibitem{kn:polya}{P\'olya,~G., 
Aufgabe 424, 
{\em Archiv der Mathematik und Physik\/} {\bf 20} (1913), pp.~271.} 

\bibitem{kn:rado}{Rado,~R., Axiomatic treatment of rank in infinite sets, 
{\em Canadian Journal of Mathematics\/} {\bf 1} (1949), pp.~337--343.} 

\bibitem{kn:rado2}{Rado,~R., 
Note on the transfinite case of hall's theorem on 
representatives, 
{\em Journal of the London mathematical Society\/} {\bf 42} (1967), 
pp.~321--324.} 
 
\bibitem{kn:redei}{R\'edei,~L., 
{\em Lacunary polynomials over finite fields\/}, 
Akademiai Kiado, 1973. \\
Also published by 
North-Holland, Amsterdam - London; 
American Elsevier, New York, 1973.} 

\bibitem{kn:riviere}{de~Rivi\`{e}re,~A., 
Question nr. 48,  
{\em L'Interm\'ediare des Math\'ematiciens\/} {\bf 1} (1894), 
pp.~19--20.}

\bibitem{kn:robertson}{Robertson,~N., P.~D.~Seymour and R.~Thomas, 
Permanents, Pfaffian orientations, and even directed circuits, 
{\em Annals of Mathematics\/} {\bf 150} (1999), pp.~929--975.} 

\bibitem{kn:ryser}{Ryser,~H.~J., 
{\em Combinatorial Mathematics\/}, 
Carus Mathematical Monographs~{\bf 14}, Published by the 
Mathematical Association of America, distributed by 
John Wiley and Sons, 
1963.}  

\bibitem{kn:sands}{Sands,~A.~D., 
The factorization of abelian groups (II), 
{\em The Quarterly Journal of Mathematics\/} {\bf 13} (1962), 
pp.~45--54.} 

\bibitem{kn:sands2}{Sands,~A.~D., 
Factoring finite abelian groups, 
{\em Journal of Algebra\/} {\bf 275} (2004), pp.~540--549.} 

\bibitem{kn:shapiro}{Shapiro,~B., M.~Shapiro and A.~Vainshtein, 
Periodic De Bruijn triangles: exact and asymptotic results, 
{\em Discrete Mathematics\/} {\bf 298} (2005), pp.~321--333.} 

\bibitem{kn:stanley}{Stanley,~R.~P., 
A survey of alternating permutations. 
ArXiv: 0912.4240v1, 2009.} 

\bibitem{kn:szabo2}{Szab\'o,S., 
Constructions related to the R\'edei property of groups, 
{\em Journal of the London Mathematical Society\/} {\bf 73} (2006), 
pp.~701--715.} 

\bibitem{kn:szabo}{Szab\'o,~S. and A.~D.~Sands, 
{\em Factoring groups into subsets\/}, 
Chapman \& Hall / CRC, Lecture Notes in Pure and Applied 
Mathematics {\bf 257}, 2009.}  

\bibitem{kn:temperley}{Temperley,~H.~N.~V. and M.~E.~Fisher, 
Dimer problem in statistical mechanics -- an exact result, 
{\em Philosophical Magazine\/} {\bf 6} (1961), pp.~1061--1063.} 

\bibitem{kn:tutte2}{Tutte,~W.~T., 
On the enumeration of planar maps, 
{\em Bulletin of the American Mathematical Society\/} {\bf 74} (1968), 
pp.~64--74.} 

\bibitem{kn:tutte}{Tutte,~W.~T. and C.~A.~B.~Smith, 
On unicursal paths in a network of degree 4, 
{\em American mathematical Monthly\/} {\bf 48} (1941), pp.~233--247.} 

\bibitem{kn:valiant}{Valiant,~L., 
The complexity of computing the permanent, 
{\em Theoretical Computer Science\/} {\bf 8} (1979), pp.~189--201.} 

\bibitem{kn:vazirani}{Vazirani,~V.~V. and M.~Yannakakis, 
Pfaffian orientations, 0-1 permanents, and even cycles 
in directed graphs, 
{\em Discrete Applied Mathematics\/} {\bf 25} (1989), pp.~179--190.} 

\bibitem{kn:viennot}{Viennot,~G., 
Permutations ayant une forme donnee, 
{\em Discrete Mathematics\/} {\bf 26} (1979), pp.~279--284.} 

\bibitem{kn:waerden}{van~der~Waerden,~B.~L., 
Ein Satz \"uber Klasseneinteilungen von endlichen Mengen, 
{\em Abhandlungen aus dem Mathematischen Seminar der Universit\"at 
Hamburg\/} {\bf 5} (1927), pp.~185--188.} 

\bibitem{kn:wilson}{Wilson,~R.~M., 
Decompositions of complete graphs into subgraphs isomorphic to a 
given graph, {\em Proceedings $5^{\mathrm{th}}$ British Combinatorial 
Conference\/}, Congressus Numerantium XV, Utilitas Mathematica (1976), 
pp.~647--659.} 

\begin{document}

\thispagestyle{empty}

\StartArtikel[Titel={\sc De Combinatoriek van De Bruijn},
          AuteurA={T. Kloks},
          AdresA={Department of Computer Science\crlf 
                  National Tsing Hua University\crlf
                  Hsinchu, Taiwan},
          EmailA={kloks@cs.nthu.edu.tw} 
	  ]
	 
\StartLeadIn
Im memoriam: N.~G.~de~Bruijn. Dit is een korte samenvatting van z"yn werk 
in de combinatoriek. 
\StopLeadIn

In dit artikel probeer ik een kort overzicht te 
geven van het werk van De Bruijn in de combinatoriek. 

\onderwerp{\underline{Representantensystemen}}

\medskip 

In 1927 publiceert Van~der~Waerden z"yn stelling over een gemeenschappel"yk 
representanten systeem van twee opsplitsingen van een eindige 
verzameling~\cite{kn:waerden}. 
Om precies te z"yn luidt z"yn stelling als volgt. 

\begin{theorem}[Van der Waerden]
Z"y $M$ een eindige verzameling. Laat $\mathcal{U}$ en 
$\mathcal{B}$ twee partities\footnote{Voor alle 
duidel"ykheid: een partitie van een verzameling 
$V$ is een verzameling deelverzamelingen van $V$ waarvan er geen  
leeg is, die elkaar paarsgew"ys niet overlappen,  en die samen elk element 
van $V$ bevatten.} z"yn van $M$ zodanig dat elk 
element van $\mathcal{U}$ en elk element van $\mathcal{B}$ 
precies $n$ elementen heeft. Laat $\mu=|\mathcal{U}|=|\mathcal{B}|$. 
Dan is er een verzameling $X \subset M$ 
van $\mu$ elementen zodanig dat voor elke 
$x \in X$ er precies \'e\'en $U \in \mathcal{U}$ is, en er precies \'e\'en 
$B \in \mathcal{B}$ is, zodat  
\[X \cap U=X \cap B=\{x\}.\] 
\end{theorem}
Met andere woorden, $X$ is een gemeenschappel"yk representantensysteem 
voor $\mathcal{U}$ en $\mathcal{B}$. 

\bigskip 

In een nawoord legt Van~der~Waerden uit dat z"yn stelling equivalent is 
met de stelling van K\H{o}nig die zegt dat elke reguliere bipartiete 
graaf een bedekking met dimeren heeft 
(zie Hoofdstuk~\ref{hoofdstuk dimeren}). 
H"y merkt op dat de stelling van K\H{o}nig naar het oneindige 
is uitgebreid 
(waarb"y de graad van de reguliere bipartiete graaf eindig bl"yft) 
in een artikel van K\H{o}nig en Valk\'o~\cite{kn:konig}. 
Uit die uitbreiding volgt ook de uitbreiding van Van~der~Waerden's 
stelling naar het oneindige. Daarb"y is dan $M$ mogel"yk oneindig, 
en elk element van $\mathcal{U}$, en elk element van 
$\mathcal{B}$, heeft precies hetzelfde eindige aantal,  
zeg $n$, elementen.

\bigskip 

Uitbreidingen van de nauw verwante 
huwel"yksstelling van Hall naar 
het oneindige volgen in 1965 door 
Rado~\cite{kn:rado2}. In het eindige geval is Hall's 
stelling equivalent met de stelling van K\"onig-Egerv\'ary. 
Deze stelling zegt dat, 
voor bipartiete grafen, het maximaal aantal elementen 
in een verzameling kanten z\'odanig dat geen 
twee kanten een punt gemeen hebben, gelijk is aan het minimale 
aantal punten wat alle kanten raakt. Hieruit volgt natuurl"yk 
K\"onig's stelling voor reguliere bipartiete grafen die we 
hierboven noemden. 

\bigskip 
 
Rado noemt ook eerdere versies van 
de uitbreiding van Hall's stelling, namel"yk  
door Hall zelf in 1948~\cite{kn:hall} en 
anderen~\cite{kn:everett,kn:halmos}. 
Deze vroege uitbreidingen (door Hall, door Everett en Whaples, 
en door Halmos en Vaughan) 
maken gebruik van Tychonoff's stelling 
(zie Hoofdstuk~\ref{hoofdstuk BE grafen}). Rado's bew"ys 
maakt g\'e\'en gebruik van Tychonoff's stelling. Bovendien 
bew"yst Rado veel meer dan alleen de uitbreiding die we 
hieronder noemen in Stelling~\ref{thm Hall inf}.   
 
\begin{theorem}[Hall]
\label{thm Hall inf}
Laat $N$ een eindige of oneindige verzameling z"yn en laat 
voor elke $k \in N$, $A_k$ een eindige verzameling z"yn. 
Stel dat aan Hall's voorwaarde is voldaan, dat wil zeggen 
dat voor elke eindige $M \subset N$ geldt 
\[|\bigcup_{k \in M} A_k| \geq |M|.\] 
Dan is er voor elke $k \in N$ een $x_k \in A_k$   
z\'odanig dat 
\[k \neq \ell \quad\text{impliceert}\quad x_k \neq x_{\ell}.\] 
\end{theorem}

\bigskip 

In 1943 versch"ynt een artikel van De~Bruijn met een 
andere uitbreiding 
van Van~der~Waerden's stelling naar het oneindige. 
Dat wil zeggen dat voor 
$M$ nu ook een oneindige verzameling mag worden gekozen. De eis dat elk 
element van $\mathcal{U}$ en van $\mathcal{B}$ hetzelfde, 
eindige aantal 
elementen bevat wordt vervangen door de volgende twee eisen. 

\begin{center}
\begin{boxedminipage}[h]{8.3cm}
\begin{enumerate}[\rm 1.]
\item
Voor geen enkele $k \in \mathbb{N}$ z"yn er 
$k-1$ elementen van $\mathcal{U}$ waarvan de 
vereniging $k$ verschillende 
elementen van $\mathcal{B}$ bevat. 

\medskip 

\noindent
Hetzelfde geldt met $\mathcal{U}$ en $\mathcal{B}$ hierboven 
verwisseld. 
\item
Elk element van $\mathcal{U}$ en van $\mathcal{B}$ 
heeft slechts met een eindig aantal elementen van respectievel"yk 
$\mathcal{B}$ 
en $\mathcal{U}$ een niet-lege doorsnede. 
\end{enumerate}
\end{boxedminipage}
\end{center}

\bigskip 

De~Bruijn's artikel~\cite{kn:bruijn24} vervolgt, na een korte 
inleiding met een zestal opmerkingen. Ik noem er drie. 
Op de eerste plaats volgt uit de eerste eis, met $k=1$,  
dat geen enkele 
deelverzameling van $B$ of $U$ leeg is. 
De tweede opmerking van De~Bruijn die ik noem stelt dat, \'als er een 
gemeenschappel"yk representantensysteem is, dat dan aan het eerste 
item is voldaan. Met andere woorden, als de tweede eis  
geldt, dan is er een gemeenschappel"yk 
representantensysteem dan en slechts dan als aan de 
eerste eis voldaan is. 
De laatste, $6^{\mathrm{de}}$ 
opmerking van De~Bruijn is, dat we mogen aannemen 
dat $\mathcal{B}$ en $\mathcal{U}$ 
(hoogstens) aftelbaar z"yn. Dit volgt uit een constructie 
van K\H{o}nig~\cite[Stelling G, op Bladz"yde 460]{kn:konig2}. 

\bigskip 

Uit dat laatste volgt dat we 
de elementen van $\mathcal{U}$ en $\mathcal{B}$ 
kunnen nummeren, wat we meteen maar doen, 
zeg 
\[U_1,U_2,\dots \quad\text{en}\quad B_1,B_2,\dots.\]  

\bigskip 

\begin{definition}
Een $S_k$-systeem is een verzameling $\{x_1,\dots,x_k\}$ van elementen 
uit $M$ zodanig dat geen element van $\mathcal{U}$ of van 
$\mathcal{B}$ er twee van bevat. 
\end{definition}

\bigskip 

\begin{lemma}
\label{lm repres}
Onder aanname van de eerste eis alleen  
geldt het volgende. 
Z"y $k \geq 0$ en z"y $\{x_1,\dots,x_k\}$ een $S_k$-systeem 
zodat 
\begin{equation}
\label{eqnrepres1}
U_i \cap B_i = \{x_i\} \quad\text{voor $\;i \in \{1,\dots,k\}$.}
\end{equation}
Dan kan men b"y elke $B_t$ met $t > k$ een $s > k$ vinden, en een 
$S_{k+1}$-systeem $\{y_1,\dots,y_{k+1}\}$ wat  
\[U_1,\dots,U_{k}, U_{s} \quad\text{en}\quad B_1,\dots,B_k,B_{t}\]  
representeert (maar niet meer pers\'e in dezelfde volgorde). 
\end{lemma}

\begin{remark}
De~Bruijn merkt op dat z"yn bew"ys grotendeels dat van 
Van~der~Waerden volgt. W"y volgen De~Bruijn's beschr"yving. 
Het construeert natuurl"yk, in de terminologie 
van Berge~\cite{kn:berge}, `een verbeterend pad'.\footnote{Zie b"yvoorbeeld 
ook~\cite[Exercise 6.7, op Bladz"yde 146 en Exercise 6.13 
op Bladz"yde 149]{kn:nienhuys}.} 
\end{remark}

\begin{bewijs}[Van Lemma~\ref{lm repres}.] 
We zeggen dat een verzameling $U_s$ verbonden is met de 
verzameling $B_t$  
als er een $\omega \geq 1$ en een ketting  
\begin{equation}
\label{eqnrepres2}
B_{t}=B_{i_0}, \;U_{i_1}, \;B_{i_1}, \;U_{i_2},\;B_{i_2}, \;\dots, 
\;B_{i_{\omega-1}},\;U_{i_{\omega}}=U_{s}     
\end{equation}
bestaan zodat 
\[B_{i_{\lambda -1}} \cap U_{i_{\lambda}} \neq \es 
\quad\text{voor $\;\lambda \in \{1,\dots,\omega\}$.}\]
We beweren dat er minstens \'e\'en $s > k$ is zodat 
$U_{s}$ verbonden is met $B_{t}$. 
Stel namel"yk dat $U_{j_1},\dots,U_{j_p}$ 
de elementen z"yn van $\mathcal{U}$ 
die verbonden met $B_t$, waarb"y 
\[j_1 < j_2 < \dots < j_p \leq k.\]   

\medskip 

\noindent
Laat dan 
$N= \cup_{z=1}^p U_{j_z}$. Beschouw $x \in B_t$ 
en laat $x \in U_m$ voor zekere $m$. Dan is $U_m$ verbonden met 
$B_t$, en dus ligt $U_m$ in $N$. Dus $B_t \subset N$.  

\medskip 

\noindent
Laat $x_z \in B_{j_z}$ voor zekere $1 \leq z \leq p$. Stel dat 
$x_z \in U_m$ voor zekere $m$. Dan is de ketting die loopt 
tot $U_{j_z}$ te verlengen met $B_{j_z}$ en $U_m$. Dus ook 
$U_m$ is met $B_t$ verbonden, en dus is $x_z \in N$. 
Met andere woorden, $B_{j_z} \subset N$, voor $1 \leq z \leq p$. 

\medskip 

\noindent
Dan is 
\[\cup_{z=1}^p B_{j_z} \cup B_t \subset \cup_{z=1}^p U_{j_z}.\] 
Dit is in str"yd met Eis~1  
en dit bew"yst de bewering. 

\medskip

\noindent
Er bestaat dus een ketting vanuit $B_t$ naar een $U_s$ met $s > k$. 
Natuurl"yk nemen we de kleinst mogel"yke $s$ waar dat gebeurt, 
zeg $\omega$, dus $i_1, \dots, i_{\omega -1}$ z"yn allemaal 
hoogstens $k$. We mogen ook aannemen dat er geen twee dezelfde elementen 
van $\mathcal{U}$ in de ketting zitten, anders kunnen we een 
deel van de ketting weglaten. Kies nu een 
element 
\[y_{\lambda} \in B_{i_{\lambda -1}} \cap U_{i_{\lambda}} \quad 
\text{voor $\;\lambda \in \{1,\dots,\omega\}$.}\] 
In het representantensysteem vervangen we 
\[x_{i_1},\dots,x_{i_{\omega}-1} \quad\text{door}\quad 
y_{i_0},y_{i_1},\dots,y_{i_{\omega}-1}.\] 
Die liggen achtereenvolgens in 
\[U_{i_1}, \dots,U_{i_{\omega -1}}, U_s \quad\text{en in}\quad 
B_{t},B_{i_1},\dots,B_{i_{\omega -1}}.\] 

\medskip

\noindent
Dit nieuwe representantensysteem voldoet aan de eis. 
\qed\end{bewijs}

\bigskip 

\begin{theorem}
\label{thm repres}
Veronderstel dat de partities $\mathcal{U}$ en $\mathcal{B}$ van 
een eindige of oneindige verzameling $M$ voldoen aan de eisen~1 en~2. 
Dan bestaat er een gemeenschappel"yk representantensysteem 
$\mathcal{U}$ en $\mathcal{B}$. 
\end{theorem}
\begin{bewijs}
We kunnen aannemen dat $M$ oneindig is en dat $\mathcal{U}$ en 
$\mathcal{B}$ aftelbaar z"yn. Neem een nummering van de 
elementen van $\mathcal{U}$ en $\mathcal{B}$, zeg 
\[U_1,\;U_2,\;\dots \quad\text{en}\quad B_1,\;B_2,\;\dots.\] 
Uit Lemma~\ref{lm repres} volgt nu het volgende. 
Als er een  
$S_k$-systeem (voor de een of andere $k$) 
$\mathcal{S}_n$ is wat onder andere $\{U_1,\dots,U_n\}$ 
en $\{B_1,\dots,B_n\}$ representeert, dan is er ook een $S_k$-systeem 
$\mathcal{S}_{n+1}$ 
wat onder andere de elementen van  
$\{U_1,\dots,U_{n+1}\}$ en van $\{B_1,\dots,B_{n+1}\}$ representeert. 

\medskip 

\noindent
We vinden een aftelbare r"y $\mathcal{S}_1,\mathcal{S}_2,\dots$ van 
$S_k$-systemen. In verschillende $\mathcal{S}_i$'s kan 
$U_1$ met verschillende $B_{\mu}$'s corresponderen 
(dat wil zeggen dat die $B_{\mu}$'s een element met 
$U_1$ gemeen hebben). 
Volgens de tweede eis is het aantal $B_{\mu}$'s wat 
$U_1$ doorsn"ydt eindig. 
Er is dus minstens \'e\'en $B_{\mu^{\prime}}$ die 
in oneindig veel $\mathcal{S}_i$'s correspondeert met $U_1$. De 
$\mathcal{S}_i$'s die deze eigenschap niet hebben laten we uit de r"y weg. 

\medskip 

\noindent
We kunnen dit proces herhalen, achtereenvolgens voor 
$U_1,B_1,U_2,B_2,\dots$. We vinden telkens een 
nieuwe correspondentie $(U_n,B_{\mu_n})$ en $(B_n,U_{\nu_n})$ 
voor zover die niet al eerder was vastgelegd. Dit leidt 
tot een \'e\'en-\'e\'en-duidige correspondentie tussen 
elementen van $\mathcal{U}$ en $\mathcal{B}$ die paarsgew"ys 
een niet-lege doorsnede hebben. 

\medskip

\noindent
Daarmee is de stelling bewezen. 
\qed\end{bewijs}

\bigskip 

\begin{remark}
Laten we de oneindige bipartiete graaf $G$ eens bekijken, 
waarvan de knopen de elementen van $\mathcal{U}$ en van $\mathcal{B}$ 
z"yn. Een paar $\{U_i,B_j\}$ vormt een kant 
in $G$ als $U_i \cap B_j \neq \es$. 
De eerste eis is dan de (tweez"ydige) voorwaarde van Hall. De tweede 
eis zegt dat elke punt een eindige graad heeft. 
Volgens Stelling~\ref{thm repres} is er dan een bedekking met dimeren, 
dat wil zeggen, er is een verzameling kanten die paarsgew"ys geen punt 
gemeen hebben en die samen alle punten van $G$ bedekken. 
Dus elke knoop kan worden uitgehuwel"ykt.  
\end{remark}

\onderwerp{\underline{De Bruijn cykels}}

\medskip 

\noindent
Wellicht is De Bruijn het meest bekend om een van z"yn eerste 
artikelen~\cite{kn:bruijn1,kn:compeau}.  
In dit artikel bew"yst De Bruijn een gissing van ir.~K.~Posthumus. 
Naar eigen zeggen kwam De~Bruijn er na het schrijven van het artikel 
achter dat het al eens eerder was bewezen, namel"yk door 
Camille Flye Sainte-Marie~\cite{kn:flye,kn:riviere,kn:bruijn2}. 

\bigskip 
 
Het originele artikel van De~Bruijn betreft 
een stelling over woorden met letters 
uit een alfabet met twee letters, 0 en 1. 
Later, samen met Van~Aardenne-Ehrenfest, heeft h"y het  
uitgebreid naar alfabetten met een willekeurig aantal 
letters~\cite{kn:aardenne}. 
Over dit artikel later meer. 

\begin{definition}
Laat $n \in \mathbb{N}$. Een $P_n$-cykel is een geordende 
cykel van $2^n$ cijfers 0 en 1 zodanig dat de $2^n$ 
r"yen van opeenvolgende nullen en enen  
allemaal verschillend zijn. 
\end{definition}
Met andere woorden, elke mogel"yke r"y nullen en enen ter lengte $n$ 
komt precies \'e\'en maal voor als een opeenvolgende deelr"y in 
de cykel. 

\bigskip 

Als voorbeeld kan men een $P_3$-cykel nemen 
\[ 00010111 \] 
waarin men de 8 mogelijke r"ytjes 
\[000 \quad 001 \quad 010 \quad 101 \quad 011 \quad 
111 \quad 110 \quad\text{en}\quad 100\] 
als deelr"ytjes kan onderkennen. 

\bigskip 

Natuurlijk is er slechts \'e\'en $P_1$-cykel, namelijk $01$, en er 
is slechts \'e\'en $P_2$-cykel, namelijk $0011$. Door wat te 
puzzelen kan men wel vinden dat er twee $P_3$-cykels 
z"yn en zestien $P_4$-cykels. Een ingenieur, ir. K.~Posthumus, ploos uit 
dat het aantal $P_5$-cykels 2048 is en op grond daarvan 
giste h"y (zie~\cite{kn:bruijn2}) de volgende stelling. 

\begin{theorem} 
\[\text{Het aantal $P_n$-cykels is} \quad 2^{2^{n-1}-n}  
\quad\text{voor alle $\;n$.}\]
\end{theorem}
\begin{bewijs}
Om dit te bew"yzen maken we gebruik van de volgende 
definitie. 

\begin{definition}[Zie ook~\cite{kn:good}]
Beschouw een gerichte graaf $(G,\Delta)$. De l"yngraaf 
$L(G,\Delta)$ is de gerichte graaf $(G^{\prime},\Delta^{\prime})$ 
gedefinieerd door 
\[G^{\prime}= \Delta \quad\text{en}\quad
\begin{cases}
\begin{minipage}{6cm} 
$\Delta^{\prime}$ bestaat uit alle paren uit $\Delta$ die 
kop-staart gerelateerd zijn. Formeel:
\end{minipage}
\\
\Delta^{\prime}=\{\;((P,Q),(R,S)) \in \Delta \times \Delta \;|\; Q=R\;\}.
\end{cases}\] 
\end{definition}
Laat nu $G_n$ de verzameling van alle woorden z"yn met $n-1$ letters. 
Voor elk element $(\epsilon_1,\dots,\epsilon_n) \in G_{n+1}$ 
nemen we  
een gerichte kant   
\[((\epsilon_1,\dots,\epsilon_{n-1}), (\epsilon_2,\dots,\epsilon_n)) 
\in \Delta_n.\] 
Een $P_n$-cykel is dan een gesloten wandeling door $(G_n,\Delta_n)$ 
die elke kant precies een maal bezoekt. Dit reduceert de vraag naar het 
aantal $P_n$-cykels tot de vraag naar het aantal Euler circuits 
in $(G_n,\Delta_n)$.  
 
\medskip 

\noindent
Merk nu op dat $(G_{n+1},\Delta_{n+1})=L(G_n,\Delta_n)$. 
Laat $N_n$ het aantal Euler circuits in $(G_n,\Delta_n)$ zijn. 
In het artikel wordt  
bewezen dat 
\[N_{n+1}=2^{2^{n-1}-1} \cdot N_n = 2^{2^{n-1}-1} \cdot 2^{2^{n-1}-n} = 
2^{2^n-n-1}.\] 
Hieruit volgt de stelling.  
\qed\end{bewijs}

\bigskip 

\begin{remark}
Dit artikel stamt uit 1946. Het is wellicht 
opmerkel"yk dat het geen enkele 
referentie bevat, iets wat men tegenwoordig b"yna niet meer ziet. 
De~Bruijn maakt dat later goed in~\cite{kn:bruijn2}. 
Het adres dat De~Bruijn geeft is het `Natuurkundig Laboratorium der 
N.V. Philips' Gloeilampenfabrieken'.  
\end{remark} 

\onderwerp{\underline{De De~Bruijn-Erd\"os stelling uit de 
incidentiemeetkunde}}

\medskip 

\noindent
Twee jaar later, in 1948, versch"ynt een artikel van 
De~Bruijn en Erd\"os.        
De stelling uit dit artikel gaat de geschiedenis in als 
de `De~Bruijn-Erd\"os stelling'.  
Er dient evenwel opgemerkt te worden 
dat er tw\'e\'e van die stellingen zijn; een in de 
grafentheorie en een in de incidentiemeetkunde. 
We beginnen met de laatste. 

\bigskip 

Voor $n \in \mathbb{N}$ en voor $m \in \mathbb{N}$ met $m >1$, 
laat 
\[U=\{1,\dots,n\},\] en laat $A_1,\dots,A_m$ deelverzamelingen 
van $U$ voorstellen. Neem aan dat elk paar elementen van $U$ bevat 
is in precies \'e\'en van de deelverzamelingen $A_i$.  
Dan geldt de volgende stelling. 

\begin{theorem}
\label{thm B-E}
Er geldt $m \geq n$ waarb"y gel"ykheid optreedt slechts dan in  
een van de volgende twee gevallen. 
\begin{enumerate}[\rm (1)]
\item Een deelverzameling bevat 
$n-1$ elementen, zeg 
\[A_1=\{1,\dots,n-1\}.\] 
De andere deelverzamelingen z"yn dan, zonder verlies 
van algemeenheid, $A_i=\{n,i-1\}$, 
voor $i=2,\dots,n$. 
\item Het getal $n$ is van de vorm $n=k(k-1)+1$. 
Alle $A_i$s hebben dan $k$ elementen en elk element komt in 
precies $k$ deelverzamelingen voor. 
\end{enumerate}
\end{theorem}

\bigskip 

\begin{corollary}
Z"y gegeven een configuratie van $n$ punten in het 
vlak die niet allemaal op \'e\'en l"yn liggen. Verbindt elk 
tweetal van die punten. Dan is het 
aantal l"ynen in dit systeem minstens $n$. In dit geval 
treedt gel"ykheid slechts dan op als 
er $n-1$ punten op een l"yn liggen.
\end{corollary}

Dit uitvloeisel van Stelling~\ref{thm B-E} is ook een 
gevolg van een stelling 
van Gallai die bekend staat als de `Sylvester-Gallai stelling'.  
In het artikel van De~Bruijn en Erd\"os wordt ook 
Gallai's elegante 
bew"ys van die stelling beschreven.

\begin{theorem}[Sylvester-Gallai stelling~\cite{kn:gallai,kn:motzkin}]
\label{thm gallai}
Stel dat $n$ punten in het vlak gegeven z"yn en dat 
niet alle punten op 
\'e\'en l"yn liggen. Dan is er een l"yn die door precies 
twee punten gaat. 
\end{theorem}
\begin{bewijs}
Stel dat de stelling onjuist is. Dan gaat elke l"yn die 
door twee punten gaat ook nog door een derde punt. Projecteer 
een van de punten, zeg $a$, in het oneindige en verbindt het 
met alle andere punten. Dan kr"ygen we een stelsel 
evenw"ydige l"ynen die elk twee of meer van de overige punten bevat. 
Neem nu een l"yn door twee punten die een zo klein mogel"yke 
hoek maakt met de evenw"ydige l"ynen. Stel dat deze l"yn de 
punten $a_1$, $a_2$ en $a_3$ bevat. De l"yn die $a_2$ met $a$ verbindt 
bevat minstens een derde punt, zeg $a_4$. Maar dan maakt een van de 
l"ynen, {\em of\/} die door $a_1$ en $a_4$,  {\em of\/} 
die door $a_3$ en $a_4$, 
een kleinere hoek 
met het stelsel evenw"ydige l"ynen. 
(Door middel van een figuur kan men zichzelf 
daarvan eenvoudig overtuigen.) 

\medskip 

\noindent
Dit is een tegenspraak die de stelling bew"yst. 
\qed\end{bewijs}

\bigskip 

Het bewijs van Stelling~\ref{thm B-E} gaat ruwweg als volgt. 

\begin{bewijs}[Van Stelling~\ref{thm B-E}.] 
We noemen de deelverzamelingen $A_1,\dots,A_m$ l"ynen 
en we noemen de elementen (of `punten'), van $U$ $a_1,\dots,a_n$. 
Laat $k_i$ het aantal l"ynen z"yn wat door het punt $a_i$ gaat 
en laat $s_j$ het aantal punten op de l"yn $A_j$ z"yn. 
Dan geldt natuurl"yk 
\begin{equation}
\label{eqn1}
\sum_{j=1}^m s_j = \sum_{i=1}^n k_i.
\end{equation}

\medskip 

\noindent
Ook geldt dat, als $a_i$ niet op de l"yn $A_j$ ligt, 
\begin{equation}
\label{eqn2}
s_j \leq k_i.
\end{equation}
Dat is zo omdat $a_i$ verbonden is door een l"yn met 
elk punt op de l"yn $A_j$, en omdat elk tweetal van die l"ynen  
verschillend zijn. 

\medskip 

\noindent
Laat $k_n$ nu de kleinste $k_i$ zijn en laat $A_1,\dots,A_{\nu}$ 
de l"ynen z"yn die door $a_n$ gaan (waarb"y $\nu=k_n$). 
We mogen aannemen dat elke 
l"yn minstens twee punten bevat, anders kunnen we die l"yn ook 
wel weglaten. Ook 
geldt dat $k_n > 1$ omdat anders alle punten op een l"yn 
zouden liggen. Neem nu, voor 
$i=1,\dots,\nu$, een punt $a_i \neq a_n$ op l"yn $A_i$. 
Dan concluderen we uit~\eqref{eqn2} dat 
\begin{eqnarray}
\label{eqn3}
&& s_2 \leq k_1 \quad s_3 \leq k_2 \quad \cdots \quad s_{\nu} \leq k_{\nu-1} 
\quad s_1 \leq k_{\nu};  
\nonumber\\
&& \quad\text{en, voor $j > \nu$,}\quad s_j \leq k_n.
\end{eqnarray}
Nu volgt uit~\eqref{eqn1} en~\eqref{eqn3}, en uit de keuze van $k_n$, 
dat $m \geq n$. 
  
\medskip 

\noindent
We analyseren nu de gevallen waarin $m=n$. Als $m=n$ dan geldt 
gelijkheid in alle ongel"ykheden van~\eqref{eqn3} gel"ykheid.  
Als $n=m$ kunnen we de punten hernummeren zodanig dat 
\[s_1=k_1 \quad \cdots \quad s_n=k_n.\]
We mogen aannemen dat 
\[k_1 \geq \dots \geq k_n > 1.\]
We onderscheiden twee gevallen. 

\begin{description}
\item[Geval $k_1 > k_2$] 
Dan is $s_1=k_1 > k_i$, $2 \leq i \leq n$. Nu volgt uit~\eqref{eqn2} 
dat alle punten op $A_1$ liggen. Het punt $a_1$ ligt dan niet op 
$A_1$ en dus geldt het eerste geval uit de stelling. 
\item[Geval $k_1=k_2$] 
Stel dat $k_j < k_1$. Dan, volgens~\eqref{eqn2}, 
ligt $a_j$ op $A_1$ en op $A_2$. De enige mogel\"ykheid is dat 
$j=n$. Omdat $k_1=\dots=k_{n-1} > k_n \geq 2$ is bevat 
elke l"yn door $a_n$, behalve \'e\'en, minstens twee andere punten. 
Er z"yn dus minstens twee l"ynen die $a_n$ niet bevatten. 
Voor beiden geldt volgens ~\eqref{eqn2} dat $s_j \leq k_n$, 
maar is in tegenspraak met $s_1=\dots=s_{n-1} > k_n$. 

\noindent
We hebben dus, behalve het hierboven genoemde geval, 
alleen het geval dat $s_i=k_j=k_1$ voor alle $1 \leq i,j \leq n$. 
Het is nu makkel"yk na te gaan dat $n=k(k-1)+1$ en dat elk paar 
l"ynen elkaar sn"ydt in precies een punt. 
\end{description}
Dit bew"yst de stelling.
\qed\end{bewijs}

\bigskip 

\begin{remark}
Een gevolg uit Stelling~\ref{thm B-E},  
die men in de grafentheorie wel vaker tegenkomt, luidt als volgt. 
Laat $\mathcal{S}$ een niet-triviale partitie z"yn van de 
kanten van $K_n$ in klieken. Dan geldt dat $|\mathcal{S}| \geq n$ en 
gelijkheid geldt dan en slechts dan als  
\begin{enumerate}[\rm (a)]
\item of een kliek $C  \in \mathcal{S}$ bevat $n-1$ knopen en 
de overige 
$n-1$ klieken z"yn kopie\"en van 
$K_2$ die elk de enige knoop bevatten die niet in $C$ ligt, 
\item of $n=k^2-k+1$ en $\mathcal{S}$ bestaat 
uit $n$ klieken 
die elk $k$ knopen bevatten. Bovendien is dan elke knoop in precies 
$k$ klieken van $\mathcal{S}$ bevat. 
\end{enumerate}
De enige kliek $K_n$ die een niet-triviale partitie van de kanten in 
$n$ driehoeken toelaat, is dus $K_7$. In dat geval  
is het Fano-vlak de incidentiemeetkunde uit Stelling~\ref{thm B-E}. 
Volgens een stelling van Wilson geldt dat, 
als $\binom{n}{2}$ deelbaar is door 
$\binom{k}{2}$, en als $n-1$ deelbaar is door $k-1$, en als $n$ 
groot genoeg is, dan is er een partitie van de kanten van $K_n$ in 
$k$-klieken~\cite{kn:wilson}.    
\end{remark}

\onderwerp{\underline{Bases voor integers}}

\medskip 

In 1950 versch"ynt een artikel~\cite{kn:bruijn20} 
wat een inspiratiebron wordt 
voor veel later werk, van zowel De~Bruijn zelf 
als van veel anderen, eg,~\cite{kn:dowell,kn:eigen}. 
Als voorbeeld noem ik een concept waar 
De~Bruijn z"yn naam aan verleent; 
de `Moser-De~Bruijn r"y'~\cite{kn:bruijn21,kn:moser}. 
Dit is de r"y getallen die de som z"yn van 
verschillende machten van vier. De r"y begint 
als volgt.  
\[0,\;1,\;4,\;5,\;16,\;17,\;20,\;21,\;64,\;65,\;\dots\] 
De getallen in deze r"y hebben allerlei interessante 
eigenschappen, en vormen het onderwerp van veel studie. 

\bigskip 

In dit Debrecen-artikel uit 1950 beantwoord De~Bruijn onder 
andere een vraag van Szele. Laten we beginnen met de definitie 
van een basis; het is natuurl"yk wat je denkt dat het is.   

\begin{definition}
Een verzameling gehele getallen 
\[\{\;b_1,\;b_2,\; \ldots\;\}\] 
is een basis voor $\mathbb{Z}$ als elke getal 
$x \in \mathbb{Z}$ op een unieke manier geschreven kan worden als 
\begin{equation}
\label{eqn13}
x=\sum_{i=1}^{\infty} \epsilon_i b_i \quad\text{waar 
elke $\epsilon_i \in \{0,1\}$ en waar 
$\sum_{i=1}^{\infty} \epsilon_i < \infty$.}
\end{equation}
\end{definition}

\bigskip 

Szele opperde het vermoeden dat elke basis 
slechts \'e\'en oneven getal heeft, slechts \'e\'en oneven veelvoud 
van twee, slechts \'e\'en oneven veelvoud van vier, enzovoort. 
Hieronder geven we De~Bruijn's bew"ys van de juistheid 
van dit vermoeden. We beginnen met een lemma. 

\bigskip 

\begin{lemma}
\label{lm szele}
Als $B=\{b_1,b_2,\dots\}$ een basis is dan is 
er \'e\'en $b_i$ oneven en alle andere z"yn even. 
\end{lemma}
\begin{bewijs}
Tenminste een $b_i$ moet oneven z"yn anders kan een oneven 
getal niet geschreven worden als combinatie van elementen uit de 
basis. Aangezien de opsomming van de elementen van $B$ willekeurig is, 
is het voldoende om te laten zien dat $b_1b_2$ even is.  

\medskip 

\noindent
Laat $V_1$ de verzameling gehele getallen z"yn waarvoor in~\eqref{eqn13} 
$\epsilon_1=0$. Laat $V_2$ de verzameling gehele getallen z"yn 
waarvoor in~\eqref{eqn13} $\epsilon_2=0$ en laat $W=V_1 \cap V_2$. 
Beschouw twee gehele getallen $x$ en $y$ met $x-y=b_1$. Veronderstel 
dat $y \in V_1$. Dan geldt  
\[y=\sum_{i=2}^{\infty} \epsilon_i(y) b_i \quad\text{wat impliceert}\quad 
x=b_1+\sum_{i=2}^{\infty} \epsilon_i(y) b_i.\] 
Dan is $x \notin V_1$  
want volgende de definitie van een basis is $x$ maar op \'e\'en 
manier te schr"yven in de vorm~\eqref{eqn13}. 
Natuurl"yk impliceert $\epsilon_1(x)=1$ dat 
$\epsilon_1(y)=0$. Met andere woorden, 
precies \'e\'en getal van $x$ en $y$ behoort tot $V_1$. 
 
\medskip 

\noindent
Als $x \in V_1$ dan z"yn, volgens de bovenstaande redenering,  
$x-2b_1$ en $x+2b_1$ ook in $V_1$. Dus $V_1$ is periodiek 
met periode $2b_1$. Zo is ook $V_2$ periodiek met periode $2b_2$. 
Dan is $W=V_1 \cap V_2$ periodiek met periode $2b_1b_2$. Laat $p$ de 
kleinste periode z"yn voor de elementen van $W$. 

\medskip 

\noindent
Voor $\lambda \in \{0,1\}$ en voor $\mu \in \{0,1\}$ z"y $W_{\lambda \mu}$ 
de verzameling getallen met 
\[\epsilon_1=\lambda \quad\text{en}\quad \epsilon_2=\mu.\] 
Merk op dat de verzamelingen $W_{\lambda \mu}$ bijectief op elkaar 
afgebeeld kunnen worden door eenvoudige verschuivingen. Z"y 
\[K=|\{\;x\;|\; x\in W \quad\text{en}\quad 0 < x \leq P \;\}|.\] 
Uit de symmetrie volgt nu dat elke set $W_{\lambda \mu}$ in 
elke periode precies  
$K$ elementen bevat. Dan volgt dat $P=4K$ en dus is $4$ een deler 
van $2b_1b_2$. Met andere woorden, $b_1b_2$ is even, 
hetgeen te bew"yzen was. 
\qed\end{bewijs}

\bigskip 

Nu bew"yzen we Szele's vermoeden. 

\begin{theorem}
\label{thm szele}
Elke basis kan geschreven worden in de vorm 
\begin{equation}
\label{eqn14}
B=\{d_1,\;2d_2,\;2^2d_3,\;\dots\;\}
\end{equation}
waarb"y de getallen $d_i$ oneven z"yn.
\end{theorem}
\begin{bewijs}
Laat $\{b_1,b_2,\dots\}$ een basis z"yn. 
Volgens lemma~\ref{lm szele} mogen we aannemen dat 
$b_1$ oneven is en dat alle andere $b_i$'s even z"yn. Nu is 
$\{b_2,b_3,\dots\}$ een basis voor de even getallen en dus is 
\[\{\;\frac{1}{2}b_2,\;\frac{1}{2}b_3,\;\frac{1}{2}b_4,\dots\;\}\] 
een basis voor de gehele getallen. Dan is dus precies 
\'e\'en van de elementen oneven. Enzovoort. 
\qed\end{bewijs}

\bigskip

\begin{remark}
Z"y $D=[d_1,d_2,\dots]$ een r"y oneven getallen. 
De r"y heet een fundament als~\eqref{eqn14} een basis is. 
Stel dat $D$ periodiek is, dat wil zeggen dat 
$d_{i+s}=d_i$ voor zekere $s > 0$ en alle $i$. 
De~Bruijn bew"yst dat dan in een eindig aantal stappen 
vastgesteld kan worden of $D$ een fundament is. 
Het artikel bevat een l"yst van alle 20 fundamenten 
$[a,b,a,b,\dots]$ met periode twee, met    
\[0 < -b < a \leq  100.\] 
In het artikel wordt Schutte  
bedankt voor het controleren van die 
l"yst fundamenten.\footnote{Professor 
H.~J.~(Hennie)~Schutte 
is in 1957 een van de stichters van 
`Die Suid-Afrikaanse Wiskundige Vereniging'.} 
In 1964 breidt De Bruijn dit 
resultaat uit~\cite{kn:bruijn21}. 
   
Aan het eind van het artikel duikt de vraag van Haj\'os op als 
(fout) vermoeden. We zien dit terug 
in Hoofdstuk~\ref{hoofdstuk factorisatie}.
\end{remark}
 
\onderwerp{\underline{De BEST stelling}}

\medskip 

\noindent
In 1951 versch"ynt het artikel van Van~Aardenne-Ehrenfest 
en De~Bruijn wat we al eerder noemden~\cite{kn:aardenne}. 
Achterin, als `note added in proof',   
wordt vermeld dat het aantal Euler circuits in een gerichte graaf 
kan worden uitgedrukt als een determinant. 
In deze notitie staat ook dat het artikel van 
Tutte en Smith ~\cite{kn:tutte} hetzelfde resultaat aankondigt. 
Om die reden gaat de stelling de geschiedenisboeken in als 
de `BEST stelling', vanwege de initialen \underline{B}ruijn, 
\underline{E}hrenfest, 
\underline{S}mith en \underline{T}utte. 
 
\bigskip 

Zij $(G,\Delta)$ een gerichte graaf. Om het aantal opspannende 
bomen, met alle kanten gericht naar een gegeven wortel, 
te vinden kan men 
gebruik maken van de stelling van Kirchoff. 

\begin{theorem}
Zij $(G,\Delta)$ een gerichte graaf met knopen  
$p_1,\dots,p_n$. Beschouw de matrix $M$ 
met 
\[M(i,j)=\begin{cases}
\text{als $i=j$:} 
\quad \text{de totale uitgraad van punt $p_i$}\\
\text{als $i \neq j$:} 
\quad -\# \;(\text{takken gericht van $p_i$ naar $p_j$}) 
\end{cases}\] 
Dan is het aantal opspannende bomen gericht naar een wortel  
$p_k$ gelijk aan de $(k,k)$-minor van $M$, dat is de 
determinant van de matrix die men krijgt door uit $M$ de 
$k^{\mathrm{de}}$ r"y en de $k^{\mathrm{de}}$ kolom weg te laten.
\end{theorem}
Voor een bew"ys z"y de lezer verwezen naar~\cite{kn:nienhuys}. 

\bigskip 

Zij nu $(G,\Delta)$ een gerichte graaf waarin voor elk 
punt de ingraad gel"yk is aan de uitgraad. Dan is $G$ een 
zogenaamde Eulerse graaf en men kan dan, 
volgens Euler's unicursal stelling, 
een wandeling, een zogenaamde 
`Euler-tour', maken door de graaf 
die elke tak precies een maal bezoekt. 
De BEST stelling telt het aantal mogel"yke Euler tours. 

\bigskip 

\begin{theorem}[De BEST stelling] 
\label{thm best}
Z"y $(G,\Delta)$ een gerichte Eulerse graaf. Laat 
de ingraad en uitgraad van een knoop $p$ gegeven zijn door 
$\sigma_p$. 
Dan is het aantal Euler-tours in $(G,\Delta)$ gelijk aan 
\[P_w(G,\Delta) \; \prod_{p \in G} \; (\sigma_p -1)! \] 
waarb"y $P_w(G,\Delta)$ het aantal opspannende bomen 
is gericht naar een wortel $w$. 
\end{theorem}

\bigskip 

\begin{remark}
Een eigenschap van Eulerse grafen is dat 
\[P_w(G,\Delta)=P_v(G,\Delta)\] 
voor elke twee knopen $v$ en $w$. 
Dus de keuze van de knoop $w$ in Stelling~\ref{thm best} is 
willekeurig.  
\end{remark}

\bigskip 

\begin{bewijs}[Van Stelling~\ref{thm best}.] 
Laat de punten van $G$ genummerd z"yn als 
$p_1,\dots,p_n$. Schr"yf $\sigma_i$ in plaats van $\sigma_{p_i}$ 
voor de ingraad van de knoop $p_i$. 

\medskip 

\noindent
Beschouw een Euler tour en start de tour 
met een gerichte 
kant $(p_1,p_2)$. Nummer de kanten die bezocht worden door de 
Euler tour en beschouw voor elke knoop $\neq p_1$ de laatste uitgang 
die genomen wordt. Het is niet moeil"yk in te 
zien dat deze kanten een opspannende boom vormen; 
de zogenaamde `laatste-uitgangs-boom.'  

\medskip 

\noindent  
We bewijzen dat elke opspannende boom 
met wortel $p_1$, en waarin elke kant 
gericht naar de wortel $p_1$, precies $\prod(\sigma_i-1)!$  
keer voor komt als laatste-uitgangs-boom. Kleur de kant 
$(p_1,p_2)$ rood. Voor elke knoop kleur de laatste uitgang 
die genomen wordt in de Euler-tour blauw. 
Dan krijgen we een blauwe boom. 

\medskip 

\noindent
Beschouw nu een opspannende boom gericht naar $p_1$. 
Kleur de kanten die 
in de boom zitten blauw. Dan heeft elke knoop precies 
\'e\'en blauw gekleurde uitgang, behalve de knoop $p_1$. 
Kleur \'e\'en uitgaande kant vanuit 
$p_1$, zeg $(p_1,p_2)$ rood. 

\medskip 

\noindent
Fixeer nu, voor elke knoop, een ordening van de uitgaande kanten, 
op een zodanige manier dat de gekleurde kant de laatste in 
de ordening is. Het aantal 
manieren om dat te doen is natuurlijk 
\[\prod_{i=1}^n (\sigma_i -1)!\]

\medskip 

\noindent
De b"ybehorende tour start met de rode kant vanuit $p_1$ 
en kiest in elke knoop de volgende uitgang in de lokale 
ordening. Het is niet moeil"yk in te zien dat dit een 
Euler tour oplevert. Stel dat we een kant vanuit een 
punt $q \neq p_1$ niet gebruikt hebben. Volg dan het blauwe pad 
vanuit $q$. Zeg dat $q^{\prime}$ de opvolger is van $q$. 
Dus $(q,q^{\prime})$ is een blauwe kant die vertrekt uit $q$. 
Aangezien er een inkomende kant in $q^{\prime}$ niet is 
gebruikt, is ook de blauwe (laatste) 
uitgaande kant uit $q^{\prime}$ niet gebruikt. We kunnen dus vanuit 
$q^{\prime}$ ons pad voortzetten. Dit pad moet eindigen in $p_1$ 
omdat de blauwe kanten een opspannende boom gericht naar $p_1$ vormen. 
Daar vinden we een tegenspraak; er is een blauwe kant die in 
$p_1$ binnenkomt niet gebruikt, en dus zijn nog niet alle 
uitgangen vanuit $p_1$ op, d.w.z., gebruikt in de Euler-tour. 

\medskip 

\noindent
Dit bew"yst de stelling. 
\qed\end{bewijs}

\bigskip 

\begin{remark}
De stelling laat zien dat het aantal Euler-tours 
in polynomiale t"yd te berekenen is. 
Dit probleem is $\#$P-compleet in ongerichte grafen~\cite{kn:brightwell}. 
\end{remark}

\onderwerp{\underline{De De~Bruijn-Erd\"os stelling uit de 
grafentheorie}}
\label{hoofdstuk BE grafen}

\medskip 

\noindent
In 1951 versch"ynt dan het tweede artikel van De~Bruijn en Erd\"os 
wat bekendheid krijgt als de `De~Bruijn - Erd\"os stelling'. 
De~Bruijn is intussen werkzaam bij de universiteit Delft. 
Het onderwerp betreft hier een knopenkleuring van oneindige grafen. 

\begin{definition}
Z"y $k \in \mathbb{N}$. 
Een $k$-kleuring van een graaf $G$ is een toekenning van  
een kleur aan iedere knoop,  
uit een verzameling van $k$ kleuren, zodanig dat  
de twee eindpunten van elke kant een verschillende kleur krijgen. 
\end{definition}

\bigskip 

\begin{theorem}
\label{thm B-E2}
Z"y $k \in \mathbb{N}$. Een graaf $G$ is $k$-kleurbaar dan en slechts 
dan als elke eindige deelgraaf van $G$ $k$-kleurbaar is. 
\end{theorem}

\bigskip 

In het korte artikel stellen de auteurs dat Stelling~\ref{thm B-E2} 
een gevolg is van de stelling van Rado~\cite{kn:rado}. 
Rado's stelling is, op zijn beurt, een eenvoudige 
toepassing van Tychonoff's stelling~\cite{kn:gottschalk}. 

\bigskip 

In het artikel van De~Bruijn en Erd\"os wordt Rado's stelling als 
volgt gepresenteerd. 

\begin{theorem}[Rado's stelling~\cite{kn:rado,kn:gottschalk}]
Z"y $M$ en $M_I$ willekeurige verzamelingen. 
Stel dat voor elke $i \in M_I$ een eindige deelverzameling 
$A_i \subset M$ gegeven is. Neem aan dat er voor elke eindige 
deelverzameling $N \subset M_I$ een keuze-functie $f_N:N \rightarrow M$ 
gegeven 
is die aan elk element $i \in N$ een element uit $A_i$ 
toekent: 
\[f_N(i) \in A_i.\] 
Dan bestaat er een keuze-functie $f:M_I \rightarrow M$, 
met $f(i) \in A_i$ als $i \in M_I$, met de volgende eigenschap. 
Voor elke eindige deelverzameling $K \subset M_I$ 
is er een eindige deelverzameling $N$ met $K \subset N \subset M_I$ 
zodanig dat 
\[f(i) = f_N(i) \quad\text{voor elk element $\;i \in K$.}\]
\end{theorem}

\bigskip 

Het artikel van De~Bruijn en Erd\"os laat vervolgens 
zien dat Stelling~\ref{thm B-E2} eenvoudig volgt uit 
Rado's stelling. 

\begin{bewijs}[Van Stelling~\ref{thm B-E2}.]
Laat $M$ de verzameling z"yn van de $k$ kleuren en laat 
$M_I$ de verzameling knopen zijn van de graaf $G$. Neem voor 
elke knoop $i \in M_I$ de  
deelverzameling $A_i$ gel"yk aan $M$. 
Met elke eindige verzameling knopen 
$N \subset M_I$ correspondeert een ge\"induceerde deelgraaf $G_N$. 
De aanname is dat al die eindige deelgrafen $G_N$ $k$-kleurbaar z"yn. 
Dat wil zeggen dat er een functie $f_N: N \rightarrow M$ is die 
$G_N$ kleurt. 

\medskip 

\noindent
De keuze-functie $f: M_I \rightarrow M$, die volgt uit Rado's stelling, 
kleurt $G$. Om dat in te zien, laat $e=\{x,y\}$ 
een kant zijn van $G$. We laten zien dat $f$ verschillende kleuren 
toekent aan $x$ en $y$. Laat $K=e=\{x,y\}$. Laat $N$ een eindige 
deelverzameling zijn van $M_I$ met 
\[K \subset N \subset M_I \quad 
\text{zodanig dat $f(z)=f_N(z)$ voor alle $z \in K$.}\] 
Aangezien $f_N$ een $k$-kleuring is van de eindige deelgraaf $G_N$ 
die de kant $e=\{x,y\}$ bevat, geldt dat $f_N(x) \neq f_N(y)$ 
en dus ook $f(x) \neq f(y)$. 

\medskip 

\noindent
Dit bew"yst de stelling. 
\qed\end{bewijs}

\bigskip 

Het artikel~\cite{kn:bruijn4} past Stelling~\ref{thm B-E2} 
vervolgens toe om de volgende stelling te bewijzen. 
We hebben de volgende definitie nodig. 

\begin{definition}
Laat $S$ een verzameling z"yn en laat voor 
elke $s \in S$ een deelverzameling $f(s)$ van $S \setminus\{s\}$ 
gegeven zijn. 
Twee elementen $b$ en $c$ van $S$ heten onafhankelijk 
als 
\[b \in S \setminus f(c) \quad\text{en}\quad c \in S \setminus f(b).\] 
Een deelverzameling $A \subset S$ heet onafhankel"yk 
als elk tweetal elementen in $A$ onafhankel"yk z"yn. 
\end{definition}

\bigskip 

Een toepassing van Stelling~\ref{thm B-E2} levert de volgende stelling op. 

\begin{theorem}
Z"y $k \in \mathbb{N}$ en neem aan dat $|f(s)| \leq k$ 
voor alle $s \in S$. Dan is $S$ de vereniging van $2k+1$ onafhankel"yke 
verzamelingen. 
\end{theorem}

\bigskip

\begin{remark}
Merk op dat de 4-kleurenstelling voor  
planaire grafen met behulp van 
Stelling~\ref{thm B-E2} uit te breiden is naar 
oneindige planaire grafen. 
Verder kan de stelling b"yvoorbeeld gebruikt worden (op soortgel"yke 
w"yze als hierboven) 
om de volgende uitbreiding van Dilworth's stelling naar oneindige 
partieel geordende verzamelingen te verkr"ygen. Een (eventueel oneindige) 
partieel geordende verzameling heeft een eindige breedte $w$ dan 
en slechts dan als de verzameling op te delen is in $w$ ketens.    
\end{remark}

\onderwerp{\underline{Wortelbomen in het platte vlak}}

\medskip 

\noindent
In een artikel uit 1964 laten Harary, Prins en Tutte zien 
dat een tweetal klassen van bomen dezelfde voortbrengende functie 
hebben en ze brengen een \'e\'en-\'e\'enduidige relatie tussen 
de twee klassen tot stand~\cite{kn:harary}. 
De beschr"yving van die relatie is nogal ingewikkeld 
en De~Bruijn en Morselt laten zien dat het ook heel 
gemakkelijk kan; en wel op drie verschillende manieren~\cite{kn:bruijn7}. 

\bigskip 

\begin{definition}
Een platte wortelboom is een boom met een wortel 
die is ingebed in het platte vlak en waarb"y een ordening 
van de uitgaande takken uit de wortel gegeven is. 
\end{definition}

\bigskip 

In de volgende stelling leiden we de voortbrengende 
functie af voor het aantal platte wortelbomen met $n$ knopen. 

\begin{theorem}
\label{thm wortelboom}
Laat de voortbrengende functie voor het aantal platte 
wortelbomen gegeven zijn door 
\[f(x)=\sum_{n=1}^{\infty} c_n x^n\]
waar $c_n$ het aantal platte wortelbomen met $n$ punten is. 
Dan voldoet $f$ aan  
\begin{equation}
\label{eqn4}
f(x)= x+f(x)^2.  
\end{equation}
Voor de co\"effici\"enten $c_n$ geldt 
\[c_n=\frac{(2n-2)!}{(n-1)! n!} \quad\text{oftewel}\quad 
c_{n+1}=\frac{1}{n+1} \binom{2n}{n}.\] 
\end{theorem}
\begin{bewijs}
Een platte wortelboom $T$ bestaat \'{o}f uit slechts een punt, 
\'{o}f er is een eerste kant, vanuit de wortel die een 
platte wortelboom draagt met, zeg $k$, punten. De rest van de 
boom $T$, inclusief de wortel, is dan een platte wortelboom met $n-k$ punten. 
Dit bew"yst dat het polynoom $f$ voldoet aan de 
functionaal vergel"yking~\eqref{eqn4} en 
hieruit volgt dat  
\begin{equation}
\label{eqn5}
f(x)=\frac{1-\sqrt{1-4x}}{2}.
\end{equation}
Uit de reeksontwikkeling van 
$\sqrt{1-4x}$ is gemakkelijk af te leiden dat 
$c_n$ het $(n-1)^{\mathrm{ste}}$ Catalan getal is. 
\qed\end{bewijs}

\bigskip 

\begin{definition}
Een binaire platte wortelboom is een 
platte wortelboom die ofwel uit slechts een punt bestaat ofwel 
een wortel van graad twee heeft en waarin dan elke andere 
knoop 
graad 1 of graad 3 heeft. 
\end{definition}

\bigskip 

\begin{theorem}
Als $g(x)$ de voortbrengende functie is van de 
binaire platte wortelbomen, dan geldt 
\begin{equation}
\label{eqn19}
g(x)=\frac{f(x^2)}{x}.
\end{equation}
\end{theorem}
\begin{bewijs}
Een binaire platte wortelboom heeft \'{o}f maar een punt, 
\'{o}f heeft een wortel met daaruit twee takken. Die twee takken 
dragen vervolgens binaire platte wortelbomen. 
Dit bew"yst dat het polynoom $g$ voldoet aan 
\begin{equation}
\label{eqn6}
g(x)=x+x g(x)^2 \quad \Rightarrow \quad 
g(x)=\frac{1-\sqrt{1-4x^2}}{2x}.
\end{equation}
Laat nu $h(x)=f(x^2) /x$. Dan volgt uit~\eqref{eqn4} 
dat 
\[h(x)=x(1+h(x)^2).\] 
Dus de oplossing voor~\eqref{eqn6} is $g(x) = h(x)$. 

\medskip 

\noindent
Men kan natuurl"yk~\eqref{eqn19} ook afleiden uit~\eqref{eqn5} 
en~\eqref{eqn6}.  
\qed\end{bewijs}

\bigskip 

\begin{corollary}
\label{corr pwb vs bb}
Het aantal platte wortelbomen met $n$ knopen is gelijk aan het 
aantal binaire platte wortelbomen met $2n-1$ knopen, dat 
wil zeggen, $n$ knopen van graad 1, $n-2$ van graad 3 
en \'e\'en wortel. 
\end{corollary}

\bigskip 

We tonen Gevolg~\ref{corr pwb vs bb}  
op een tweede manier aan. 
 
We beginnen weer met de platte wortelbomen. Stel dat een luis 
van links naar rechts loopt, over de boom. B"y elke tak roept 
hij, of zij, of hij naar boven loopt of naar beneden. We krijgen 
dan van de luis een zogenaamde `up-down code' door. 

\bigskip 

B"yvoorbeeld, b"y een platte wortelboom die geen enkele 
kant heeft roept de luis niets. Een platte wortelboom 
met twee knopen heeft ook slechts \'e\'en mogelijke code, namel"yk 
`UD'. Er z"yn twee mogel"yke platte wortelbomen met drie 
knopen. De twee mogel"yke codes z"yn `UUDD' en `UDUD'. 

\bigskip 

We kunnen de algemene vorm van een UD-code weergeven met de 
Backus-Naur Form grammatica. Dat wil zeggen, 
een UD-code ziet er in het algemeen uit als volgt. 
\[<UD-code> ::= | U <UD-code> D <UD-code> \] 
Het is of een leeg woord, of het begint met een `U', gevolgd door 
een willekeurige UD-code en dan een `D', gevolgd door weer een 
willekeurige UD-code. 
  
\bigskip 

\begin{remark}
Laat $d_i$ het aantal U's z"yn die komen na de $(i-1)^{\mathrm{ste}}$ 
D en voor de $i^{\mathrm{de}}$ D. De r"y parti\"ele sommen 
voldoen dan aan 
\[d_1 \geq 1 \quad d_1+d_2 \geq 2 \quad 
d_1+d_2+d_3 \geq 3 \quad d_1+\dots+d_n=n.\] 
Het aantal mogel"yke r"ytjes van dit soort parti"ele sommen 
wordt gegeven door de Catalan getallen. 
\end{remark}

\bigskip 

De platte binaire wortelbomen kunnen we 
coderen we met een zogenaamde `knoop-eind code'. 
De luis loopt weer over de platte boom van links naar rechts. 
B"y elke knoop roept de luis of het een 
interne knoop is, zeg `K', of een blad, zeg `E'. 
Echter, h"y roept dat all\'e\'en 
maar t"ydens de \'e\'erste keer dat h"y de knoop ziet! 

\bigskip 

B"yvoorbeeld, als de boom slechts uit een wortel bestaat 
roept de luis alleen `E'. 
Als de boom bestaat uit een wortel met twee bladeren, roept 
hij `KEE'. Er twee mogel"yke 
platte binaire wortelbomen met v"yf knopen. De codes z"yn `KKEEE' en 
`KEKEE'. 

\bigskip 

De KE-code eindigt alt"yd met een `E'. Als we die weglaten 
kr"ygen we de zogenaamde afgekorte KE-code, ofwel de `aKE-code'. 

\bigskip 

Als we de KE-code met de Backus-Naur grammatica weergeven krijgen 
we de volgende formule. 
\[<KE-code> ::= E | K <KE-code> <KE-code>\] 
De vorm voor de afgekorte KE-code is dus 
\begin{eqnarray*}
<aKE> &::=& | K < KE-code> <aKE-code>\\
& ::=& | K <aKE>E<aKE>
\end{eqnarray*}
en we zien dus dat de grammatica's voor de aKE-code en de 
UD-code identiek z"yn.      

\bigskip 

\begin{remark}
Het artikel van De~Bruijn en Morselt laat ook nog op een derde  
manier de overeenkomst tussen de twee soorten bomen zien, namelijk door 
middel van een simpel plaatje dat, in feite, alles zegt.    
We verw"yzen verder ook naar~\cite{kn:nienhuys} 
voor een aardige beschr"yving van het een en ander. 

Tutte liet het trouwens niet op zich zitten; die kwam een 
paar jaar later met de telling van `platte grafen'~\cite{kn:tutte2}.
Daarin noemt h"y wel nog het artikel van Harary, Prins en zichzelf, 
maar niet dat van De~Bruijn en Morselt.
\end{remark}

\bigskip 

\begin{remark}
In 1972 versch"ynt nog een artikel van De~Bruijn over platte wortelbomen.  
Het kr"ygt een behoorl"yke invloed in de loop der t"yd 
vanwege 
toepassingen, onder andere in de informatica.  
Dit keer is het geschreven samen met Knuth en Rice en gaat het over 
de gemiddelde hoogte van platte wortelbomen met $n$ knopen. 
Na een fikse rekenpart"y is het resultaat de volgende stelling. 

\begin{theorem}
Gemiddeld genomen over alle platte wortelbomen met $n$ knopen is 
de hoogte  
\[ \sqrt{\pi \cdot n} - \frac{1}{2} + O(n^{-1/2}\log n) 
\quad (n \rightarrow \infty). \] 
\end{theorem}

\end{remark}

\onderwerp{\underline{Permutaties met een vaste vorm}}

\medskip 

Permutaties van een bepaalde vorm worden al 
meer dan 120 jaar 
bestudeerd~\cite{kn:andre,kn:andre2,kn:macmahon,kn:shapiro}. 
Over de alternerende permutaties is natuurl"yk het meest bekend. 
Een permutatie $[a_1,\dots,a_n]$ heet alternerend als 
\begin{equation}
\label{eqn17}
a_1 > a_2 < a_3 > a_4 < \dots 
\end{equation}
en als de tekens de andere kant op staan heet de permutatie 
omgekeerd alternerend. Als $E_n$ het aantal alternerende 
permutaties is dan geldt voor de voortbrengende 
functie~\cite{kn:andre,kn:stanley}  
\begin{multline}
\sum_{i=0}^{\infty} E_n \frac{x^n}{n!} = \sec(x)+\tan(x)= 
\frac{1}{\cos(x)}+\tan(x)=\\
= 1+x+\frac{x^2}{2!}+2\frac{x^3}{3!}+5\frac{x^4}{4!}+
16\frac{x^5}{5!}+61\frac{x^6}{6!}+\dots
\end{multline}
Ook geldt de recurrente betrekking 
\[E_{n+1}=\sum_{\substack{1 \leq j \leq n\\ \text{$j$ odd}}} 
\binom{n}{j} E_j E_{n-j} \quad \text{voor $\;n \geq 1$.}\] 
Een asymptotische benadering is~\cite{kn:andre3,kn:stanley}   
\[\frac{E_n}{n!}= \frac{4}{\pi}\left(\frac{2}{\pi}\right)^n + 
O\left(\left(\frac{2}{3\pi}\right)^n\right).\] 
  
\bigskip 

MacMahon stelt in 1908 voor om de algemene 
vorm van een permutatie te bestuderen. 

\begin{definition}
Z"y $n \in \mathbb{N}$. Laat  
\[Q=(q_1,\;\dots,\;q_{n-1})\] 
een vector z"yn met alle componenten $q_i \in \{-1,+1\}$. 
Een permutatie 
\[\pi=[\pi_1,\;\pi_2,\;\dots,\;\pi_n]\] 
van $\{1,\dots,n\}$ heeft de vorm $Q$ als 
\[q_i(\pi_{i+1}-\pi_i) > 0 \quad 
\text{voor alle $\;i \in \{1,\dots,n-1\}$.}\] 
\end{definition}

Het aantal permutaties van de vorm $Q$ geven we aan met 
$\psi(Q)$. 

\bigskip 

Niven leidt in~\cite{kn:niven} een formule 
af voor $\psi(Q)$ in de vorm van een determinant (de formule 
staat ook 
op~\cite[Pagina~190]{kn:macmahon2}). 
De~Bruijn 
leidt recurrente betrekkingen af voor het 
aantal permutaties van de vorm $Q$ die eindigen met een bepaald 
c"yfer~\cite{kn:bruijn23}. 

\begin{definition} 
Laat $Q$ een vorm z"yn voor permutaties van de graad $n$, 
dat wil zeggen, we beschouwen een deelverzameling van $S_n$, de groep 
permutaties van 
\[\{1,\dots,n\}.\] 
Laat $1 \leq j \leq n$. 
Definieer $\theta(Q;j)$ als het aantal permutaties 
$[a_1,\dots,a_n]$ van de vorm $Q$ met $a_n=j$. 
\end{definition}

Er geldt natuurl"yk  
\[\sum_{j=1}^n \theta(Q;j)=\psi(Q).\] 
  
\bigskip 

\begin{lemma}
\label{lm ups and downs}
Laat $Q=\{q_1,\dots,q_{n-1})$ een vorm z"yn en definieer 
\[(Q,1)= (q_1,\dots,q_{n-1},1).\] 
Dan geldt $\theta((Q,1);1)=0$ en, voor $1 < j \leq n+1$ 
\begin{equation}
\label{eqn15}
\theta((Q,1);j)=\sum_{1 \leq h < j} \theta(Q;h).
\end{equation}
\end{lemma}
\begin{bewijs}
Voor een permutatie  
\[\alpha=[a_1,\;\dots,\;a_{n+1}]\]
definieer 
$\ell(\alpha)=a_{n+1}$. 
Definieer $\pi(\alpha)$ als de permutatie van $\{1,\dots,n\}$ 
die men kr"ygt door 1 af te trekken van 
elke $a_i$ met $1 \leq i \leq n$ die groter is dan $a_{n+1}$. 
B"yvoorbeeld,
\[\alpha=[2,4,5,1,6,3] \quad\Rightarrow\quad 
\pi(\alpha)=[2,3,4,1,5].\] 
\medskip 

\noindent
Als $\alpha$ de vorm $(Q,1)$ heeft dan 
heeft $\pi(\alpha)$ de vorm 
$Q$. Ook geldt dan dat 
\[1 \leq \ell(\pi(\alpha)) < \ell(\alpha).\] 

\medskip 

\noindent
Z"y nu $\beta$ een permutatie van $\{1,\dots,n\}$ 
van de vorm $Q$ met $\ell(\beta) < j$. Er is precies \'e\'en 
permutatie $\alpha$ van $\{1,\dots,n+1\}$ van de vorm $(Q,1)$ 
met 
$\ell(\alpha)=j$ en $\pi(\alpha)=\beta$. Tel namel"yk 1 op 
b"y elk element wat groter is dan $j-1$, en zet een $j$ op het eind. 

\medskip 

\noindent
Dit bew"yst~\eqref{eqn15}.
\qed\end{bewijs}

\bigskip 

\begin{lemma}
\label{lm2 ups and downs}
Laat $Q=(q_1,\dots,q_{n-1})$ een vorm zijn en definieer 
\[(Q,-1)=(q_1,\;\ldots,\;q_{n-1},\;-1).\] 
Dan is $\theta((Q,-1);n+1)=0$ en voor $1\leq j \leq n$ 
\begin{equation}
\label{eqn16}
\theta((Q,-1);j)=\sum_{j-1 < h \leq n} \theta(Q;j).
\end{equation}
\end{lemma}
\begin{bewijs}
Merk op dat $\theta(Q;k)=\theta(-Q;n+1-k)$ en dus  
\[\theta((Q,-1);j)=\theta((-Q,1);n+2-j).\] 
De Formule~\eqref{eqn16} volgt nu uit~\eqref{eqn15}.
\qed\end{bewijs}

\bigskip 

De Lemmas~\ref{lm ups and downs} en~\ref{lm2 ups and downs} 
leiden tot een eenvoudig algoritme om $\psi(Q)$ uit te rekenen 
voor een gegeven vorm $Q$. Dit algoritme 
vergt $O(n^2)$ stappen. De~Bruijn bescrh"yft dit 
algoritme in~\cite{kn:bruijn23}. 

\medskip

Niven bewees in z"yn artikel ook dat voor alle $n \in \mathbb{N}$ 
\begin{equation}
\label{eqn18}
\psi(Q) < \psi(Q_0)=\psi(-Q_0) 
\end{equation}
tenz"y $Q=Q_0$ of $Q=-Q_0$ is,
waarb"y $Q_0$ de vorm is van een alternerende 
permutatie~\eqref{eqn17}, dat wil zeggen,  
\[Q_0=(1,-1,1,-1,\dots,(-1)^n).\] 
Een algoritme voor het bepalen van $\psi(Q_0)$, `het Euler getal', 
was eerder al eens beschreven door Entringer~\cite{kn:entringer}. 
De~Bruijn komt door een analyse van z"yn algoritme (wat veel l"ykt 
op Entringer's algoritme) tot een 
eenvoudig bew"ys voor dezelfde ongel"ykheid~\eqref{eqn18}. 

\medskip 

Viennot doet het allemaal nog eens over in 1979~\cite{kn:viennot}. 
H"y beweert dat z"yn algoritme het eenvoudigst is, maar het is 
van hetzelfde laken een pak. Viennot bew"yst ook Niven's 
ongel"ykheid~\eqref{eqn18} nog een keer.   

\medskip 

\begin{remark}
De~Bruijn geeft, 
in een appendix, de ALGOL-60 code voor z"yn algoritme. 
Dit algoritme is getest op de Electrologica-X8 machine 
van de THE (een machine die gemaakt werd in Nederland). 
Het programma is gedraaid voor alle 
$2^n$ vormen van permutaties met $n \leq 9$. 
ALGOL is een programmeertaal die 
ontwikkeld werd door een internationaal 
team, waaronder exponenten Bakus, Naur en D"ykstra. 
Het operating system voor de 
Electrologica-X8 was het `THE multiprogramming 
system', ontwikkeld door een team 
geleid door D"ykstra. 
De enige programmertaal die ondersteund werd 
door D"ykstra's operating system was ALGOL-60. 
\end{remark}
 
\onderwerp{\underline{Het bedekken van grafen met dimeren}}
\label{hoofdstuk dimeren}

\medskip 

\noindent
\begin{definition}
Z"y $(G,\Gamma)$ een graaf. Een bedekking van $G$ met 
dimeren is een deelverzameling $\Gamma^{\prime}$ van kanten 
zodanig dat elke knoop een eindpunt is van precies \'e\'en 
kant in $\Gamma^{\prime}$. 
\end{definition}

\bigskip 

Natuurl"yk kan niet elke graaf bedekt worden 
met dimeren. Een noodzakel"yke voorwaarde is dat het aantal 
knopen even is maar dat is in het algemeen niet voldoende. 

\bigskip 

Z"y $(G,\Gamma)$ een graaf met $n$ knopen. Voor het gemak, stel dat 
\[G=\{1,\dots,n\}.\] 
Een cykelbedekking van $G$ 
is een verzameling van cykels die onderling knoop-disjunct z"yn en die 
samen alle knopen van $G$ bevatten. We laten hier ook 
cykels van lengte twee toe; 
dit z"yn cykels die bestaan uit \'e\'en kant.  

\bigskip 

In een cykelbedekking 
heeft elke knoop $i$ precies \'e\'en `opvolger', zeg ${\sigma(i)}$, 
waarb"y \[\sigma \in S_n\] dan een permutatie is van $\{1,\dots,n\}$. 
Om het begrip  van een `opvolger' te vangen 
vatten we elke kant van $G$ op als een vereniging van twee 
gerichte kanten; voor elke kant 
$\{i,j\} \in \Gamma$ z"yn er dan twee 
gerichte kanten; $(i,j)$ en $(j,i)$. Geef elke gerichte kant 
$(i,j)$ een gewicht, zeg $a_{i,j}$, en definieer 
het gewicht van een cykelbedekking      
met b"ybehorende permutatie $\sigma$, als  
\[\text{gewicht van $\sigma$}= \prod_{i=1}^n a_{i,\sigma(i)}.\] 

\bigskip 

Beschouw de verbindingsmatrix $A$ van $G$ en vervang elke $1$ door 
een gewicht $a_{i,j}$.  Dan is de permanent van $A$ 
\[\text{permanent van $A$}= 
\sum_{\sigma \in S_n} \prod_{i=1}^n a_{i,\sigma(i)}.\]
Met andere woorden, de permanent van $A$ telt de gewichten van alle 
cykelbedekkingen van $G$ op. 

\bigskip

\begin{remark}
Voor een bipartiete graaf is het aantal bedekkingen met 
dimeren gel"yk aan de permanent van de bipartiete verbindingsmatrix. 
\end{remark}

\begin{remark}
Het l"ykt erop dat het 
berekenen van de permanent van een matrix veel moeil"yker 
is dan het berekenen van de determinant. Het berekenen van de 
permanent is namel"yk $\#$P-volledig~\cite{kn:valiant}. 
Dit houdt 
onder meer in dat  
uit het bestaan van een polynomiaal algoritme $P=NP$ volgt. 
Via een formule van Ryser kan de permanent van een $n \times n$ matrix 
berekend worden in $O(2^n n)$ tijd~\cite{kn:ryser} en er 
is op het moment niet veel beters~\cite{kn:bjorklund}. 
\end{remark}

Kasteleyn vond in 1961--1963 dat het aantal bedekkingen 
van een planaire graaf met dimeren berekend kan worden met behulp 
van de Pfaffiaan van een bepaalde antisymmetrische 
matrix~\cite{kn:kasteleyn3,kn:kasteleyn}. 
In het artikel~\cite{kn:kasteleyn2} laat h"y dat eerst zien voor 
rechthoekige 
roosters. Onafhankel"yk van Kasteleyn laten ook Temperley en Fisher 
zien dat het aantal bedekkingen van rechthoekige roosters 
met dimeren uitgedrukt kan worden in Pfaffianen~\cite{kn:temperley}.   
De reden is dat een planaire graaf op een bepaalde manier, 
zogenaamd `Pfaffiaans',  
geori\"enteerd kan worden. 
Een dergel"yke Pfaffiaanse ori"entatie brengt het probleem 
terug tot het uitrekenen van een determinant. 
We besch"yven de manier waarop De~Bruijn 
dit laat zien in~\cite{kn:bruijn19}.  
     
\bigskip 
 
\begin{definition}
Z"y $(G,\Gamma)$ een graaf en laat $G=\{1,\dots,n\}$. 
Neem aan dat elke 
kant $\{i,j\}$ een gewicht $b_{i,j}$ heeft,  
waarb"y we aannemen dat 
\[b_{i,j}=b_{j,i}.\] 
Een circuit is een deelverzameling van kanten, met ten minste  
\'e\'en element, die cyclisch geordend kunnen 
worden als   
\begin{equation}
\label{eqn7}
[\;\{i_1,i_2\},\;\{i_2,i_3\},\;\dots,\;\{i_{k-1},i_k\},\;\{i_k,i_1\}\;], 
\end{equation}
en waarb"y $i_1,\dots,i_k$ allemaal 
verschillend z"yn.\footnote{Een `circuit' is hier 
dus equivalent met een `cykel',  
zoals we dat hier boven beschreven. 
Het gebruik van beide begrippen zoals hier is ongebruikel"yk.}
De lengte van het circuit is $k$.  
Het gewicht van het circuit is  
\begin{equation}
\label{eqn8}
\prod_{\ell=1}^{k} b_{i_{\ell},i_{\ell+1}} 
\quad\text{waarin $\;i_{k+1}=i_1$.}
\end{equation}
\end{definition}
Merk op dat een circuit van lengte twee bestaat uit \'e\'en 
kant 
$\{i,j\}$ en dat het gewicht van zo'n cykel $b_{i,j}^2$ is. 

\bigskip 
 
Een circuitbedekking is een verzameling circuits 
waarvan de knopen een partitie van $G$ vormen.   
Als alle circuits in een circuitbedekking lengte 
twee hebben is het dus een bedekking met dimeren. 
Het gewicht van een 
circuitbedekking is het product van de gewichten van de circuits 
erin. Ingeval het een bedekking met dimeren 
betreft noemen we de wortel uit het gewicht 
het wortelgewicht. 

\bigskip 

\begin{definition}
Een ori\"entatie van een graaf $(G,\Gamma)$ is een antisymmetrische 
$n \times n$ matrix $\epsilon$ 
met elementen $\epsilon_{i,j}$ die voldoen aan  
\[\epsilon_{i,j}= 
\begin{cases}
0 & :\quad \text{als $\{i,j\} \notin \Gamma$, en}\\
-1 \quad\text{of}\quad +1 & :\quad \begin{minipage}[t]{4cm}
als $\{i,j\} \in \Gamma$, en wel zodanig dat 
$\epsilon_{i,j}=-\epsilon_{j,i}$. \end{minipage}  
\end{cases}\]     
\end{definition}
   
\bigskip 

Het teken van een circuit van even lengte zoals 
gegeven in~\eqref{eqn7} is 
\begin{equation}
\label{eqn12}
-\epsilon_{i_1,i_2} \cdot \epsilon_{i_2,i_3} \cdot 
\dots \cdot \epsilon_{i_k,i_1}
\end{equation}
Merk op dat 
dit teken niet afhangt van de gekozen orientatie omdat het 
circuit even is. Merk ook op dat een circuit van lengte twee 
alt"yd een positief teken heeft. 
 
\bigskip 

\begin{definition}
Een ori\"entatie is Pfaffiaans als in elke circuitbedekking met 
even circuits alle circuits een positief teken hebben. 
\end{definition} 

\bigskip 

\begin{theorem}
Z"y $\epsilon$ een Pfaffiaanse ori\"entatie van een graaf $(G,\Gamma)$. 
Laat 
\begin{equation}
\label{eqn9}
S=\sum_{\beta} \rho(\beta) 
\end{equation}
waarin gesommeerd wordt over alle 
bedekkingen met dimeren $\beta$ en waarin $\rho(\beta)$ 
het wortelgewicht is van zo'n bedekking. 
Dan geldt 
\begin{equation}
\label{eqn10} 
S^2= det(A)
\end{equation}
waarby $A$ de antisymmetrische matrix is met elementen 
\[a_{i,j}=\epsilon_{i,j} \cdot b_{i,j}.\] 
\end{theorem}
\begin{bewijs}
Natuurlijk is 
\[S^2=\sum_{(\beta_1,\beta_2)} \rho(\beta_1) \rho(\beta_2)\] 
waarb"y gesommeerd wordt over alle geordende paren 
$(\beta_1,\beta_2)$ van bedekkingen met dimeren. 
Als $\beta_1$ en $\beta_2$ bedekkingen z"yn met 
dimeren dan  
is de vereniging $\beta_1\cup \beta_2$ een 
circuitbedekking met even circuits. 
De kanten in de circuits van lengte $> 2$ z"yn 
om en om 
kanten van $\beta_1$ en kanten van $\beta_2$. 

\medskip 

\noindent
Omgekeerd, z"y $\gamma$ een circuitbedekking met even circuits. 
Z"y $\nu(\gamma)$ het aantal circuits in $\gamma$ 
van lengte $>2$. Dan z"yn er $2^{\nu(\gamma)}$ manieren om 
paren $\beta_1$ en $\beta_2$ te kiezen. Als een circuit 
van $\gamma$ lengte twee heeft dan kiezen we die kant natuurl"yk 
zowel in $\beta_1$ als in $\beta_2$. 

\medskip 

\noindent
Hiermee is bewezen dat 
\[S^2=\sum_{\gamma} 2^{\nu(\gamma)} w(\gamma)\] 
waarin gesommeerd wordt over alle circuitbedekkingen 
met even circuits en waarin $w(\gamma)$ het 
gewicht is van $\gamma$. 

\medskip 

\noindent
Volgens de Leibniz formule is de determinant van $A$ gelijk aan 
\begin{equation}
\label{eqn11}
\text{determinant $A$}= \sum_{\pi \in S_n} sgn(\pi) \cdot 
a_{1,\pi(1)} \dots a_{n \pi(n)}
\end{equation}
waar gesommeerd wordt over alle permutaties 
$\pi \in S_n$. De factor $sgn(\pi)$ is $1$ of $-1$ al naargelang 
het aantal even cykels in 
$\pi$ even of oneven is. Als $\{i,\pi(i)\} \notin \Gamma$ 
dan kunnen we $\pi$ weglaten uit de sommatie, want $b_{i,\pi(i)}=0$. 
In het b"yzonder laten we alle permutaties $\pi$ buiten beschouwing 
waarvoor $\pi(i)=i$ voor de een of andere knoop $i$. 

\medskip 

\noindent
Elke permutatie splitst uit naar 
een of meer cykels. Voor de permutaties die we beschouwen vormen de 
cykels circuits in $G$ voorzien van een bepaalde omlooprichting.  
Voor de circuits 
van lengte $> 2$ zijn er twee mogel"yke omlooprichtingen. 
Als $\gamma$ een circuitbedekking is met even circuits 
dan z"yn er $2^{\nu(\gamma)}$ omlooprichtingen en die 
corresponderen met verschillende permutaties $\pi$. 

\medskip 

\noindent
Als we de omlooprichting  
van een oneven circuit omkeren, dan wisselt het teken in~\eqref{eqn11} 
want de matrix $A$ is antisymetrisch. Dus circuitbedekkingen 
die een oneven circuit bevatten doen elkaar paarsgewijs teniet. 

\medskip 

\noindent
We beweren nu dat de termen in~\eqref{eqn11},  
die horen b"y een en dezelfde 
even circuitbedekking $\gamma$,    
allemaal gel"yk z"yn aan 
\[w(\gamma)=b_{1,\pi(1)} \dots b_{n,\pi(n)}.\] 
Namel"yk, $sgn(\pi)$ is het produkt 
van een aantal $-1$ termen wat gel"yk is aan het aantal 
(even) circuits. Omdat $\gamma$ Pfaffiaans is 
geldt dat elk circuit in een even circuitbedekking 
$+1$ is, dus het produkt van de $\epsilon$'s in elk 
circuit in $\gamma$ is, volgens~\eqref{eqn12}, $-1$.     
\qed\end{bewijs}

\bigskip 

In de rest van het artikel~\cite{kn:bruijn19} 
schetst De~Bruijn het bew"ys 
van Kasteleyn~\cite{kn:kasteleyn3,kn:kasteleyn}  
dat alle planaire grafen een Pfaffiaanse ori\"entatie hebben. 
Het komt erop neer dat er een ori\"entatie 
nodig is die voor elke even cykel $C$ waarvoor  
$G \setminus C$ een bedekking heeft met dimeren, 
voor elke omlooprichting 
van $C$, het aantal kanten van $C$ waarvan de ori\"entatie 
samenvalt met de omlooprichting oneven is. 

\bigskip 

\begin{remark}
Beschouw een planaire graaf $(G,\Gamma)$ 
die ligt ingebed in het platte vlak. 
Er is een ori\"entatie van $G$ z\'odanig dat de begrenzing 
van elk 
eindig gebied een oneven aantal kanten heeft wat met de klok meedraait. 
Dit is een Pfaffiaanse ori\"entatie. 

\medskip 

Om een dergel"yke ori\"entatie te vinden kan men, ruwweg, als 
volgt te werk gaan. Neem een opspannende boom $T$ en ori\"enteer 
de kanten van $T$ willekeurig. Neem een tweede opspannende boom 
$T^{\prime}$ op de duale van $G$ door  
twee aangrenzende gebieden 
te verbinden door een kant in $T^{\prime}$ als de grens 
g\'e\'en kant is in $T$. 
Kies het buitengebied als wortel van $T^{\prime}$. 

\medskip 

Start nu met de bladeren van $T^{\prime}$ en werk naar de wortel 
toe. Onderweg, ori\"enteer de 
kanten in $G$ die corresponderen met de 
kanten van $T^{\prime}$ z\'odanig dat elke gebied een oneven 
totaal aantal kanten krijgt wat met de klok meedraait. 
\end{remark}

\bigskip 

Volgens Kasteleyn heeft 
een ori\"entatie zoals hierboven beschreven de eigenschap   
dat voor elke cykel de pariteit 
van het aantal kanten wat met de klok meedraait  
tegengesteld is aan de pariteit van het aantal knopen 
wat is ingesloten door de cykel. Als de cykel onderdeel uitmaakt 
van een circuitbedekking met even circuits dan is het aantal punten 
wat ingesloten is natuurlijk 
even en dat is dus volgens Kasteleyn precies het geval als de 
cykel een oneven aantal 
kanten heeft wat met de klok meedraait.  

\bigskip 

De~Bruijn geeft aan dat het bew"ys voor het 
bestaan van een dergel"yke ori\"entatie terug te brengen is tot dat voor   
veelhoeken die opgedeeld z"yn in driehoeken. Voor deze laatste 
klasse van grafen is, volgens De~Bruijn,   
gemakkel"yk met inductie naar het aantal kanten   
te bew"yzen dat ze een 
Pfaffiaanse ori\"entatie hebben. De inductiestap 
is het weglaten van 
\'e\'en kant van de buitenste rand van  
de veelhoek. 

\bigskip 

Als de graaf een `duale' heeft dan volgt Kasteleyn's stelling  
ook uit het volgende lemma van Little, toegepast op 
de 
duale graaf~\cite[Lemma 1]{kn:little}.  

\begin{lemma}
Laat $(G,\Gamma)$ een eindige, samenhangende graaf z"yn en 
laat $v \in G$. Dan is er een ori\"entatie van $G$ 
zodanig dat elke knoop behalve $v$ een oneven uitgraad heeft. 
\end{lemma} 
\begin{bewijs}
Laat $\epsilon$ een ori\"entatie z"yn zodanig dat het aantal 
knopen in $G\setminus \{v\}$ wat een even uitgraad heeft minimaal is. 
Stel dat er een knoop $u \neq v$ is met een even uitgraad. 
Aangezien $G$ samenhangend is, is er een pad $P$ van $u$ naar $v$ 
in $(G,\Gamma)$. Keer in elke kant $\{i,j\}$ in $P$ het teken 
$\epsilon_{i,j}$ om. 
Dat verandert de pariteit van de uitgraad all\'e\'en 
voor de knoop $u \in G\setminus \{v\}$. 
\qed\end{bewijs}

\bigskip 

\begin{remark}
In~\cite{kn:little2} geeft Little een bewijs 
voor de stelling dat een graaf die geen deelgraaf heeft 
die homeomorf is met $K_{3,3}$,  
een Pfaffiaanse ori\"entatie heeft. 
\end{remark}

\bigskip 

Laat $A$ een vierkante matrix z"yn met elementen $0$ en $1$. 
Laat $B$ een matrix z"yn die verkregen is 
uit $A$ door enkele enen in $A$ te vervangen door $-1$. 
Dan heet $B$ een P\'olya-matrix voor $A$ als  
\[\text{permanent van $A$} = \text{determinant van $B$.} 
\quad\text{Zie~\cite{kn:polya}.}\] 
Vazirani en Yannakakis bew"yzen de volgende stelling~\cite{kn:vazirani}. 

\begin{theorem}
Een bipartiete 
graaf $G$ heeft een Pfaffiaanse ori\"entatie dan en slechts 
dan als de bipartiete verbindingsmatrix een P\'olya-matrix heeft. 
\end{theorem}

\bigskip 

\begin{remark}
In~\cite{kn:little3} geeft Little 
een nodige en voldoende voorwaarde 
opdat een bipartiete graaf een Pfaffiaanse ori\"entatie heeft. 
Robertson, Seymour en Thomas 
geven een structurele karakterisering, wat leidt tot een 
polynomiaal herkenningsalgoritme~\cite{kn:robertson}. 
\end{remark}

\onderwerp{\underline{Factorisaties van eindige groepen}}
\label{hoofdstuk factorisatie}

\medskip 

\noindent
De~Bruijn geeft aan, in z"yn inleiding op het dictaat combinatoriek, 
dat de groepentheorie niet tot de combinatoriek 
behoort~\cite{kn:nienhuys}. 
(Volgens van~Lint behoort de groepentheorie, alsook de 
combinatoriek, tot 
de discrete wiskunde~\cite{kn:lint}.)
De grenzen van de combinatoriek z"yn vaag, en  
wellicht kan combinatoriek nog het best gedefinieerd worden als 
dat wat De~Bruijn onderwees en wat 
beschreven wordt in z"yn college dictaat~\cite{kn:nienhuys}. 

Maar alle gekheid op een stokje! Ik wil toch 
enkele belangrijke resultaten in de groepentheorie 
van De~Bruijn noemen. 

\bigskip 

In 1953, De Bruijn is dan werkzaam b"y de Universiteit van Amsterdam, 
versch"ynen er twee artikelen over de factorisatie van 
eindige groepen. Deze artikelen hebben tot op 
de dag van vandaag een behoorl"yke invloed 
(zie b"yvoorbeeld~\cite{kn:amin,kn:dinitz,kn:sands,kn:sands2,kn:szabo}) 
en ik wil iets ervan kort beschr"yven. 

\bigskip 

In dit hoofdstuk is een groep synoniem met een eindige abelse groep. 

\bigskip 

Als $G$ een direct product is van cyclische groepen van 
orde $2^3$, $2$ en 
$5$ dan zeggen we dat $G$ van het type $\{2^3,2,5\}$ is. 
Laat $G$ een groep z"yn en laat $A$ en $B$ deelverzamelingen 
van $G$ z"yn. We schr"yven $G=AB$ als elk element $g$ van $G$ 
op een unieke manier te schr"yven is als $g=ab$ 
met $a \in A$ en $b \in B$. In dat geval zeggen we dat  
$G=AB$ een ontbinding is van $G$ in factoren $A$ en $B$. 
Als $H_1$ en $H_2$ subgroepen z"yn van $G$ dan betekent $G=H_1H_2$ 
dat $G$ het directe product is van $H_1$ en $H_2$. 

\bigskip 

Als $A$ een deelverzameling is van een groep $G$ en als 
$g \in G$ is dan schrijven 
we 
\[Ag=\{\;ag\;|\; a \in A\;\}.\] 
Een deelverzameling $A$ van een groep $G$ heet periodiek 
als er een element $g \in G$, $g \neq e$, is zodat $Ag=A$. 
Als $A$ en $B$ twee deelverzamelingen z"yn van 
$G$ dan schr"yven we 
\[AB=\{\;ab\;|\; a \in A \quad b \in B\;\} \quad 
\text{all\'e\'en als $|AB|=|A|\cdot|B|$}.\] 
Dat wil zeggen, we gebruiken de notatie 
$AB$ all\'e\'en als elk element van $AB$ op een unieke 
manier te schrijven is als $ab$ ($a \in A$ en $b \in B$). 

\bigskip 

Een conjecture van Minkowski werd bewezen door Haj\'os 
met behulp van de volgende stelling. 

\begin{theorem}[Haj\'os' stelling~\cite{kn:hajos}]
Stel dat $G$ een eindige abelse groep is 
en stel dat $G$ een factorisatie $G=A_1 \cdots A_m$ heeft waarb"y elke 
factor $A_i$ van 
de volgende vorm is.  
\[A_i=\{\;e,\;a_i,\;a_i^2,\dots,\;a_i^{m_i}\;\}.\]
Dan is tenminste \'e\'en van de $A_i$'s periodiek. 
\end{theorem}
      
\bigskip 

De stelling van Haj\'os leidde hem (en ook De~Bruijn) 
tot de volgende vraag: 
\begin{quote}
Als $G$ een groep is, met een ontbinding in factoren $G=AB$, 
is dan alt"yd tenminste een van de factoren periodiek? 
\end{quote}
Het antwoord (van Haj\'os) 
is `nee' en dit leidde tot de vraag welke groepen 
die `Haj\'os eigenschap' hebben. 

\begin{definition}
Een eindige abelse groep $G$ heeft de Haj\'os eigenschap als 
voor elke ontbinding $G=AB$ van $G$ in factoren 
$A$ en $B$, tenminste een factor $A$ of $B$  
periodiek is. 
\end{definition}

\bigskip 

Het antwoord is intussen bekend; 
Sands geeft de volgende karakterisering~\cite{kn:sands} (met belangr"yke 
b"ydragen van De~Bruijn en R\'edei~\cite{kn:redei}). 

\begin{theorem}
Een eindige abelse groep $G$  
heeft de Haj\'os eigenschap als $G$ isomorf is met een 
subgroep van een groep die van een van de volgende types is, 
waarb"y $p<q<r<s$ priemgetallen z"yn en $n \in \mathbb{N}$: 
$\{p^n,q\}$, $\{p^2,q^2\}$, $\{p^2,q,r\}$, 
$\{p,q,r,s\}$, $\{p,p\}$, $\{p,3,3\}$, $\{3^2,3\}$, $\{p^3,2,2\}$, 
$\{p^2,2,2,2\}$, $\{p,2^2,2\}$, $\{p,2,2,2,2\}$, $\{p,q,2,2\}$, 
$\{2^n,2\}$ en $\{2^2,2^2\}$.   
\end{theorem}

\bigskip 

Voor de cyclische groepen die de Haj\'os eigenschap 
hebben werd een karakterisering gevonden in~\cite{kn:sands2} 
(met behulp van constructies door De~Bruijn). 

\bigskip 

\begin{remark}
Een nauw verwant probleem is de classificatie van groepen met de 
R\'edei eigenschap. Een groep $G$ heeft de R\'edei eigenschap als 
voor elke ontbinding $G=AB$, met $e \in A \cap B$, 
tenminste een van $A$ en $B$ 
bevat is in een echte ondergroep van $G$ 
(dus $<A> \neq G$ of $<B> \neq G$, 
waarb"y $<A>$ de kleinste ondergroep van $G$ is die $A$ bevat). 
Voor het gemak zegt men dat $G=\{e\}$ \'o\'ok de R\'edei eigenschap 
heeft. Er z"yn intussen wat resultaten bekend~\cite{kn:dinitz,kn:szabo2}. 
\end{remark}
  
\onderwerp{\underline{P\'olya's tel-theorie}}

\medskip 

\noindent
In 1964 verschijnt een artikel van De~Bruijn waarmee h"y 
zich, met allure, op internationaal gebied  
manisfesteert als leermeester in de combinatoriek~\cite{kn:bruijn8}. 
Om het een `overzichtsartikel' te noemen 
zou dit artikel jammerl"yk tekort doen! 
Het artikel bevat 
niet alleen een kristalheldere uiteenzetting van P\'olya's 
theorie maar het zit bovendien bomvol voorbeelden, toepassingen, 
uitbreidingen enz. 
Dat het in de smaak valt moge trouwens bl"yken uit het feit dat h"y 
vaak gevraagd wordt het  
nog eens over te doen, zowel in het Duits, als in het 
Engels, en Nederlands~\cite{kn:bruijn10,kn:bruijn11,kn:bruijn15}. 

\medskip 

Dat De~Bruijn een groot bewonderaar is van George P\'olya steekt h"y 
niet onder stoelen of banken! In diverse 
artikelen~\cite{kn:bruijn9,kn:bruijn26,kn:bruijn12,kn:bruijn14,kn:bruijn13}, 
waarin h"y telkens P\'olya's stelling op niet-triviale manier uitbreidt, 
noemt h"y P\'olya's stelling steevast `de fundamentele stelling', 
en ik geloof niet dat h"y het {\em ooit\/} `de Redfield-P\'olya stelling' 
genoemd heeft. 
Z"yn `liefste boek' is, sinds de t"yd dat h"y b"y Philips 
werkt, P\'olya and Szeg\"o's Aufgaben und Lehrs\"atze aus der Analysis 
(Springer, 1925)~\cite{kn:bruijn16}. 
Z"yn twee belangr"yke artikelen over Penrose's  
vlakverdelingen draagt h"y op aan P\'olya~\cite{kn:bruijn17}.  

\medskip 

Het zou te ver gaan om al het werk van De~Bruijn met betrekking 
tot P\'olya's tel-techniek hier uiteen te zetten. 
Veel kan men vinden in Nienhuys' beschr"yving van De~Bruijn's  
college {\em combinatoriek\/}~\cite{kn:nienhuys}. 
Ik beperk me hier tot een korte uiteenzetting van De~Bruijn's 
uitleg van P\'olya's stelling zoals h"y dat doet in~\cite{kn:bruijn8} 
en in~\cite{kn:nienhuys}. 

\bigskip 

Beschouw een kubus. Er z"yn natuurl"yk $2^6$ manieren om de 
z"yvlakken te kleuren met twee kleuren, zeg wit en zwart. 
Als we de kubus draaien in de ruimte dan z"yn veel van 
die kleuringen 
niet meer van elkaar te onderscheiden. 
Wat z"yn de verschillende zwart-wit kleuringen? 

\begin{enumerate}[\rm 1)]
\item Alle z"yvlakken wit. 
\item \'E\'en vlak zwart. 
\item Twee zwarte z"yvlakken die \'e\'en ribbe delen.  
\item Twee zwarte z"yvlakken, tegenover elkaar. 
\item Drie zwarte z"yvlakken, die b"y elkaar komen in een hoekpunt. 
\item Drie zwarte z"yvlakken, in een U-vorm. 
\item Twee tegenover elkaar liggende witte z"yvlakken. 
\item Twee witte z"yvlakken die een ribbe delen. 
\item Precies \'e\'en wit z"yvlak. 
\item Alle z"yvlakken zwart. 
\end{enumerate}

Als we twee kleuringen van de kubus hetzelfde noemen als 
de ene in de andere over gaat door een draaiing van de kubus 
in de ruimte, dan z"yn er dus tien verschillende kleuringen.

\bigskip 

Het hangt duidelijk af van de groep van draaiingen die de z"yvlakken 
van de kubus 
in elkaar doet overgaan. 
We zeggen dat die groep van draaiingen {\em werkt\/} op 
de verzameling $D$ van z"yvlakken van de kubus. Dat betekent 
dat er een homomorfisme is, de `werking van de groep',  
\begin{multline}
\label{eqnpolya5}
\pi: G \rightarrow S_D \quad \text{zodanig dat} \\
\forall_{g_1 \in G}\forall_{g_2 \in G}\; 
\pi(g_1g_2^{-1})=\pi g_1(\pi g_2)^{-1}.
\end{multline}
    
\bigskip 

Belangr"yk is nu het {\em type\/} van de 
draaiing. 

\begin{definition}
Elke permutatie $\pi \in S_n$ splitst $\{1,\dots,n\}$ 
op in een collectie cykels. Het type van de permutatie $\pi$ 
is de r"y 
\[(b_1(\pi),\;b_2(\pi),\; \dots),\] 
waar $b_i(\pi)$ 
het aantal cyckels is van lengte $i$. 
\end{definition}

\bigskip 

Als voorbeeld geven we hieronder de verschillende draaiingen 
van de kubus met het type van de werking op de z"yvlakken.
Er z"yn natuurl"yk $6 \times 4=24$ verschillende draaiingen. 

\begin{enumerate}[\rm (a)]
\item Alle rotaties over 90$^\circ$, rond een as 
die door de middens van twee tegenover elkaar 
liggende z"yvlakken gaat, 
l"yken op elkaar. 
Het type is $(2,0,0,1,0,\ldots)$ en het aantal draaiingen 
van dit type is 6. 
\item We kunnen de kubus draaien 
over 180$^\circ$, rond eenzelfde as. Het type 
van zo'n draaiing is $(2,2,0,\ldots)$. Er z"yn 3 van deze draaiingen. 
\item 
We kunnen de kubus 120$^\circ$ draaien, rond een as door twee 
diagonaal tegenover elkaar 
liggende hoekpunten. Het type is $(0,0,2,\ldots)$ en er z"yn er 8 van.
\item Neem een as door 
het midden van twee tegenover elkaar liggende ribben, en draai 
de kubus 180$^\circ$. Het type is $(0,3,0,\ldots)$ en er z"yn 6 
van zulke draaiingen. 
\item Tenslotte moeten we de identieke afbeelding 
niet vergeten. Er is er een, en het type is natuurl"yk 
$(6,0,0,\ldots)$. 
\end{enumerate}
Ter controle kan men de aantallen optellen; in totaal hebben we 
\[24=6+3+8+6+1 \quad\text{draaiingen.}\]  

\bigskip 

\begin{definition}
Voor een groep $G$ met een werking $\pi$ is de 
cykelindex de veelterm (met variabelen $x_1$, $x_2$, enz.)   
\begin{equation}
\label{eqnpolya1}
P_{G,\pi}(x_1,x_2,\dots)=\frac{1}{|G|} \sum_{g \in G} 
x_1^{b_1(\pi(g))}x_2^{b_2(\pi(g))} \cdots 
\end{equation}
waarb"y $(b_1(\pi(g)),b_2(\pi(g)),\dots)$ het type is van de 
permutatie $\pi(g)$. 
\end{definition}

\bigskip 

Dus, voor de groep draaiingen van de kubus $G$, 
waarb"y $\pi$ de werking is op de z"yvlakken $D$,  
is de cykelindex, volgens bovenstaande tabel:  
\[P_{G,\pi}(x_1,x_2,\dots)=\frac{1}{24} 
(x_1^6+3x_1^2x_2^2+6x_1^2x_4+6x_2^3+8x_3^2).\] 

\bigskip 

We kunnen de groep draaiingen ook laten werken op de 
\underline{hoekpunten} van de kubus. We kr"ygen dan een andere 
cykelindex, namel"yk 
\[P_{G,\pi^{\prime}}(x_1,x_2,\dots)=
\frac{1}{24}(x_1^8+9x_2^4+6x_4^2+8x_1^2x_3^2).\] 
De uitwerking op de \underline{ribben} 
geeft de volgende cykelindex. 
\[P_{G,\pi^{\prime \prime}}(x_1,x_2,\dots)=
\frac{1}{24}(x_1^{12}+3x_2^6+6x_4^3+6x_1^2x_2^5+8x_3^4).\] 
De rechtgeaarde 
ingenieur pluist natuurl"yk allerlei regeltjes uit die helpen 
b"y de controle van zoiets. B"yvoorbeeld is de som van de co\"effici\"enten 
alt"yd 24 (dat wil zeggen, in het algemeen, $|G|$). 
 
\bigskip 

Als we willen weten wat het aantal manieren is om de 
z"yvlakken van de kubus met twee kleuren te kleuren, 
waarb"y twee kleuringen hetzelfde z"yn als ze in elkaar over gaan 
door een draaiing van de kubus, 
dan geeft de stelling van P\'olya het antwoord, namel"yk,    
\begin{multline}
\label{eqnpolya10}
P_{G,\pi}(2,2,\dots)=\frac{1}{24}(2^6+3\cdot2^22^2+6\cdot2^22+6\cdot2^3+
8\cdot2^2)\\
\frac{1}{24}(64+48+48+48+32)=\frac{240}{24}=10
\end{multline}
precies wat we hierboven ook vonden, maar dan met de hand. 

\bigskip 

We zouden ook geinteresseerd kunnen z"yn in het aantal 
(niet-equivalente) kleuringen van de z"yvlakken met twee kleuren, 
zodanig 
dat er 4 zwarte en 2 witte z"yvlakken zijn. 
Dat vinden we door voor $x_i=z^i+w^i$ in te vullen. 
We vinden 
\begin{multline}
\label{eqnpolya9}
\frac{1}{24}((z+w)^6+3(z+w)^2(z^2+w^2)^2+6(z+w)^2(z^4+w^4)\\
+6(z^2+w^2)^3+8(z^3+w^3)^2).
\end{multline}
De co\"effici\"ent van $z^4w^2$ (4 zwarte en 2 witte z"yvlakken) is 
\[\frac{1}{24}(15+9+6+18+0)=2.\] 
Met andere woorden, volgens de stelling van P\'olya kunnen we
de z"yvlakken van de kubus op 2 manieren met 4 zwarte en 2 witte 
z"yvlakken kleuren. En dat klopt!     

\bigskip 

Om de stelling van P\'olya te bew"yzen hebben we een lemma 
nodig wat vroeger het `Lemma van Burnside' heette. Door toedoen 
van De~Bruijn~\cite{kn:bruijn13} 
heet dat nu het `Lemma van Cauchy-Frobenius'. 

\bigskip 

\begin{lemma}[Cauchy-Frobenius]
Z"y $G$ een groep en stel dat $\pi$ de werking is van $G$ 
op een verzameling objecten $D$. Dus $\pi: G \rightarrow S_{D}$ 
is een homomorfisme zodanig dat voor alle $g_1 \in G$ en $g_2 \in G$   
\begin{equation}
\label{eqnpolya2}
\pi(g_1g_2)=\pi(g_1) \pi(g_2) \quad \text{en}\quad \pi(e)=I_D,
\end{equation}
waarb"y $\pi$ het eenheidselement $e \in G$ afbeeldt op 
de identieke permutatie $I_D$. Merk trouwens op dat~\eqref{eqnpolya2} 
hetzelfde is als~\eqref{eqnpolya5}.  

Definieer een equivalentierelatie $\sim$ op $D$ als volgt. 
\begin{equation}
\label{eqnpolya3}
d_1 \sim d_2 \quad\Leftrightarrow\quad 
\exists_{g \in G}\; \pi(g)d_1=d_2.
\end{equation}
dan is het aantal equivalentieklassen, $|D/G|$, gelijk aan 
\begin{eqnarray}
\label{eqnpolya4}
|D/G|&=&\frac{1}{|G|} \sum_{g \in G} \psi(g), \quad\text{waarb"y} 
\nonumber\\
\psi(g)&=&|\{\;d \in D\;|\; \pi(g)d=d\;\}|.
\end{eqnarray}
\end{lemma}

\medskip 

Dus $\psi(g)=b_1(\pi(g))$, het aantal cykels in $\pi(g)$ ter lengte 1. 

\medskip 

\begin{bewijs}
Tel het aantal paren $(g,s)$ met 
\[g \in G \quad\text{en}\quad s \in D \quad\text{en}\quad 
\pi(g) s=s\] 
op twee manieren. Voor elke $g \in G$ z"yn er $\psi(g)$ 
elementen $s \in D$ die voldoen. Dus het aantal paren 
is gel"yk aan 
\[\nu=\sum_{g \in G} \psi(g).\] 
Laat, voor $s \in D$, $\eta(s)$ het aantal $g \in G$ met 
$\pi(g)s=s$ z"yn. Dan geldt dus 
\begin{equation}
\label{eqnpolya6}
\sum_{s \in D} \eta(s)=\sum_{g \in G} \psi(g).
\end{equation}

\medskip 

Voor een $s \in D$ vormen de elementen $g \in G$ waarvoor 
$\pi(g)s=s$ een ondergroep van $G$, zeg $G_s$. 
Het aantal elementen van die ondergroep is $|G_s|=\eta(s)$. 

\medskip 

Als nu $s^{\prime} \sim s$, dan is 
het aantal elementen $g \in G$ met $\pi(g)s=s^{\prime}$ 
gel"yk aan $\eta(s)$, 
want 
\[\pi(h)s^{\prime}=s \quad\Rightarrow\quad 
(\;\pi(g)s=s^{\prime} \quad\Leftrightarrow\quad    
hg \in G_s \;).\] 

\medskip 

\noindent
We kunnen $G$ dus opsplitsen in deelverzamelingen met  
$|G_s|$ elementen z\'odanig dat elke deelverzameling 
precies \'e\'en element 
van de equivalentieklasse van $s$ bevat. 
De equivalentieklasse van $s$ bevat dus $\frac{|G|}{\eta(s)}$ 
elementen. Met andere woorden 
\[\eta(s)=\frac{|G|}{\# \;\{\;s^{\prime}\;|\; s^{\prime} \sim s\;\}}.\] 
De som over $\eta(s)$ over alle $s$ in \'e\'en equivalentieklasse 
is $|G|$. Dus de som van $\eta(s)$ over alle $s \in D$ 
is gel"yk aan $|G|$ maal het aantal equivalentieklassen. 

\medskip 

\noindent
Het lemma volgt nu uit~\eqref{eqnpolya6}. 
\qed\end{bewijs}

\bigskip 

Op precies dezelfde manier kunnen we ook een gewogen versie 
van Cauchy-Frobenius afleiden. Stel dat elk element $d \in D$ 
een gewicht $\omega(d)$ heeft. Een eis daarb"y is wel dat het 
gewicht constant is op de equivalentieklassen, dat wil zeggen, 
dat voor elk tweetal elementen 
\[d_1 \sim d_2 \quad\Rightarrow\quad \omega(d_1)=\omega(d_2).\]  
We kunnen dan het gewicht van een equivalentieklasse definieren 
als het gewicht van een element uit die klasse. 

\begin{lemma}[Cauchy-Frobenius met gewichten]
\label{lm gewogen CF}
De som van de gewichten van de equivalentieklassen $D/G$
is $\frac{1}{|G|} \sum_{g \in G} \Psi(g)$, 
waarb"y 
\[\Psi(g)=\sum_{d \in D} \omega(d) \cdot \nu(\pi(g)d=d).\] 
\end{lemma}
 
Hier gebruiken we de indicatorfunctie 
\[\nu(\text{bewering})= 
\begin{cases} 
0 & \text{als de bewering onwaar is}\\
1 & \text{als de bewering waar is.}
\end{cases}\] 

\bigskip 

\begin{definition}
Een kleuring is een functie $f: D \rightarrow R$, 
waarb"y $R$ een verzameling kleuren is, b"yvoorbeeld 
\[R=\{\;\text{wit}, \;\text{zwart}\;\},\] 
en waarb"y $D$ de verzameling objecten is die we 
willen kleuren. In ons voorbeeld is $D$ de verzameling 
z"yvlakken van de kubus.
\end{definition}

\bigskip 

Zoals al eerder opgemerkt z"yn we niet zozeer ge\"interesseerd 
in het aantal kleuringen (dat z"yn er $2^6$), maar in het aantal 
{\em kleurpatronen\/}. Twee kleuringen z"yn equivalent als 
de ene over gaat in de andere door een draaiing van de kubus, 
en een kleurpatroon is een equivalentieklasse onder die 
relatie. 

\bigskip 

Z"y $G$ weer de groep draaiingen van de kubus. 
De werking van $G$ op de kleuringen is een homomorfisme 
gedefinieerd als volgt. 
\begin{multline}
\label{eqnpolya7}
\sigma: G \rightarrow S_{R^D} \quad\text{gedefinieerd door}\\
\sigma(g)f=f \circ \pi(g^{-1}).
\end{multline}

\bigskip 

Laten we eerst even controleren of dit inderdaad een 
werking is volgens~\eqref{eqnpolya2} (of~\eqref{eqnpolya5}). 
\begin{eqnarray*}
\sigma(g_1g_2)f &=& f \circ \pi((g_1g_2)^{-1}) \\
&=& f \circ \pi(g_2^{-1}g_1^{-1}) \\
&=& f \circ \pi(g_2^{-1}) \circ \pi(g_1^{-1})\\
&=& \sigma(g_1) (f \circ \pi(g_2^{-1}) \\
&=& (\sigma(g_1) \sigma(g_2)) f.
\end{eqnarray*}

\bigskip 

Geef de kleuren een `gewicht', b"yvoorbeeld $w$ en $z$. 
Het gewicht van een kleur $r \in R$ geven we aan met 
$\omega(r)$. 

\begin{definition}
Z"y $f:D \rightarrow R$ een kleuring. Het gewicht van $f$ 
is 
\[\Omega(f)=\prod_{d \in D} \omega(f(d)).\]
\end{definition}

B"yvoorbeeld, als we de z"yvlakken van de kubus kleuren met 
wit en zwart, dan heeft elke kleuring met 4 zwarte en twee 
witte z"yvlakken hetzelfde gewicht, namel"yk $z^4w^2$. 
  
\bigskip 

Als we de kubus draaien, verandert het gewicht niet, met andere 
woorden, we kunnen het gewicht van een kleurpatroon $F$ defini\"eren 
als 
\[\Omega(F)=\Omega(f) \quad\text{voor $\;f \in F$.}\] 

\bigskip 

\begin{theorem}[P\'olya's Fundamentele Stelling]
\label{thm polya}
De som van de gewichten van de kleurpatronen is 
\begin{equation}
\label{eqnpolya8}
P_{G,\pi}(\sum_{r \in R} \omega(r), \sum_{r \in R} \omega(r)^2, 
\sum_{r \in R} \omega(r)^3,\dots).
\end{equation}
\end{theorem}

\bigskip 

In ons voorbeeld kr"ygen we dus de veelterm~\eqref{eqnpolya9}.  
Hieruit kunnen we dan b"yvoorbeeld het aantal kleurpatronen 
met twee witte en vier zwarte z"yvlakken aflezen. 
Als we alle gewichten gewoon 1 nemen, krijgen we 
\[\sum_{r \in R} \omega(r)^i=|R| \quad\text{voor alle $i$.}\] 
Het totaal aantal kleurpatronen met twee kleuren is 
dus $P_{G,\pi}(2,2,\dots)=10$, zoals in~\eqref{eqnpolya10}. 

\bigskip 

\begin{bewijs}[Van Stelling~\ref{thm polya}.]
Volgens Cauchy-Frobenius (met gewichten) geldt 
\begin{eqnarray}
\label{eqnpolya11}
\sum_{F \in \{\;\text{kleurpatronen}\;\}} \Omega(F)= \frac{1}{|G|} 
\sum_{g \in G} \Psi(g), \nonumber\\
\text{waarb"y} \; \Psi(g)=\sum_{f \in R^D} \Omega(f) \cdot \nu(\sigma(f)=f).
\end{eqnarray}
Laat $g \in G$. Dan splitst de werking $\pi(g)$ van $g$, $D$ in 
cykels $D_1$, $D_2$, enz. De eis dat $f \circ \pi(g^{-1})=f$ 
betekent dat alle cykels $D_i$ egaal gekleurd moeten worden. 
Dat houdt in 
\[\Psi(g)=\sum_{\substack{f \in R^D\\ \sigma(f)=f}} \Omega(f)=
\prod_{\text{cykels $D_i$}} \; \sum_{r \in R} \omega(r)^{|D_i|}.\]  
Als we nu de cykels b"y elkaar vegen die dezelfde lengte hebben, 
vinden we 
\begin{multline}
\Psi(g)=\left(\sum_{r \in R} \omega(r) \right)^
{\# \text{1-cykels in $\pi(g)$}} \cdot \\ 
\cdot \left(\sum_{r \in R} \omega(r)^2 \right)^
{\# \text{2-cykels in $\pi(g)$}} \cdot \dots 
\end{multline}
En dit is precies~\eqref{eqnpolya8}. 
\qed\end{bewijs}

\bigskip 

\begin{remark}
De~Bruijn publiceerde veel uitbreidingen 
op P\'olya's stelling. B"yvoorbeeld zou men 
ge"interesseerd kunnen z"yn in het aantal kleurpatronen 
waarb"y men ook kleuren mag verwisselen. 
Dit l"ykt 
minder gekunsteld als de kleuren veel op elkaar l"yken, 
zoals purper en violet; dit z"yn kleuren die veel mensen 
wel kunnen onderscheiden maar toch vaak verwisselen. 
We kr"ygen dan een soort 
`kleur-designs', 
in plaats van gewone kleurpatronen.  

\bigskip 

In b"yvoorbeeld~\cite{kn:bruijn8} en~\cite{kn:nienhuys} kan 
men nog veel 
meer uitbreidingen en voorbeelden vinden. 
\end{remark}

\bigskip 

\onderwerp{\underline{Afsluiting}}

\medskip 

Ik heb in dit artikel geprobeerd een kort overzicht te geven 
van De~Bruijn's b"ydragen tot de combinatoriek. Voor m"yzelf was 
dit een b"yzonder aangename reis door de t"yd. De artikelen van 
De~Bruijn hebben niets 
verloren aan helderheid; men leest ze alsof ze gisteren door De~Bruijn 
geschreven werden. Dat was m"yn ervaring ook, zelfs in versterkte mate, 
b"y het vertalen in het Engels van Nienhuys' college-aantekeningen. 
Ik waande me meer dan eens terug in de t"yd, luisterend naar De~Bruijn; 
luisterend naar {\em z"yn\/} combinatoriek. 

\bigskip 

Ik moest voor dit korte overzicht natuurl"yk een keuze maken. 
Er z"yn nog veel andere verhalen die het lezen zeker waard z"yn. 
Om nog \'e\'en ander, b"yzonder leuk, voorbeeld te noemen is er het 
blokkendoos verhaal uit 1969~\cite{kn:bruijn25,kn:nienhuys}. 
Daarin gaat De~Bruijn in op de vraag van  
z"yn zoon, die er op zeven-jarige leeft"yd achter komt dat h"y 
z"yn $6 \time 6 \times 6 \times 6$ doos niet kan vullen met blokken 
van $1 \times 2 \times 4$. 

\bigskip 

Ook de artikelen van De~Bruijn over vloerbedekkingen 
\`{a} la Penrose 
laat ik onbesproken. Over dit onderwerp, alsook over 
De~Bruijn's asymptotiek en z"yn werk aan AUTOMATH, versch"ynen 
binnekort andere samenvattingen. 

\bigskip 

Ik heb geprobeerd een indruk te geven van z"yn meest bekende 
artikelen in de combinatoriek, 
en ik laat het hierb"y; met een diepe buiging voor 
een der grootste Nederlandse wiskundigen \'en leermeesters.             

\completepublications

\end{document}